\documentclass[11pt]{article}
\usepackage{graphicx,epsfig}
\usepackage{amssymb}
\usepackage{amsmath}
\usepackage{theorem}
\usepackage{epsfig}
\usepackage{verbatim}
\usepackage{graphicx}
\usepackage{amsfonts}

\usepackage{color}

\newtheorem{thm}{Theorem}[section]
\newtheorem{cor}[thm]{Corollary}
\newtheorem{prop}[thm]{Proposition}

\newtheorem{rem}[thm]{Remark}

\newtheorem{lemma}[thm]{Lemma}

\newcommand{\R}{\Bbb{R}}

\newcommand{\D}{\displaystyle}

\newcommand{\di}{{\rm div}\thinspace}
\newcommand{\curl}{{\rm curl}\thinspace}
\newcommand{\supp}{{\rm supp}\thinspace}
\newcommand{\grad}{\nabla}

\newcommand{\dt}{\frac{d}{dt}}

\newcommand{\dxt}{\partial_{x_3}}
\newcommand{\dpt}{\partial_t}

\normalsize

\newcommand{\la}{\Lambda}

\newcommand{\al}{\alpha}
\newcommand{\alu}{\alpha_1}
\newcommand{\ald}{\alpha_2}

\newcommand{\dau}{\partial_{\alpha_1}}
\newcommand{\dad}{\partial_{\alpha_2}}
\newcommand{\dbu}{\partial_{\beta_1}}
\newcommand{\dbd}{\partial_{\beta_2}}
\newcommand{\dai}{\partial_{\alpha_i}}
\newcommand{\daj}{\partial_{\alpha_j}}

\numberwithin{equation}{section}

\begin{document}


\author{Antonio C\'ordoba, Diego C\'ordoba and Francisco Gancedo}
\title{Porous media: the Muskat problem in 3D}
\date{May 11, 2010}

\maketitle

\begin{abstract}
The Muskat problem involves  filtration of two incompressible
fluids throughout a porous medium. In this paper we shall discuss in
3-D the relevance of the Rayleigh-Taylor condition, and the topology
 of the initial interface, in order to prove its local existence in Sobolev spaces.
\end{abstract}

\section{Introduction}

The Muskat problem (see ref. \cite{Muskat} and \cite{bear}) involves  filtration of two incompressible fluids throughout  a porous medium characterized by a positive constant $\kappa$ quantifying its porosity and permeability.
  The two fluids, having respectively velocity fields $v^j$, $j=1,2$,  occupy disjoint regions   $D^j$ ($D^2 = \mathbb{R}^3 - D^1$) with a common boundary (interface) given by the surface $S = \partial D^j$. Naturally those domains  change with time,  $D^j = D^j(t)$, as it does its  interface $S=S(t)$. We shall denote by $p^j$ ($j=1,2$) the corresponding pressure and we will assume also that the dynamical viscosities $\mu^j$ and the densities $\rho^j$ are constants such that $\mu^1 \neq \mu^2$, $\rho^1 \neq \rho^2$.

  The conservation of mass law in this setting is given by the equation $\nabla\cdot v=0$ (in the distribution sense) where $v= v^1\chi_{D^1} + v^2\chi_{D^2}$.

  The momentum equation was obtained experimentally by Darcy \cite{Darcy, bear} and reads as follows \begin{equation*}\label{darcylaw}
\frac{\mu^j}{\kappa} v^j =-\nabla p^j - (0,0,\,\rho^j g),\quad j=1,2,
\end{equation*} where $g$ is the acceleration due to gravity.

One can find in the literature several attempts to derive Darcy law from
Navier-Stokes (see \cite{Tartar} and \cite{Sanchez}) throughout the process of
homogenization under the hypothesis of a periodic, or almost periodic, porosity. In any case the presence of the porous medium justify the elimination of the inertial terms in the motion, leaving friction (viscosity) and gravity as the only relevant forces, to which one has to add  pressure as it appears in the formulation
of Darcy's law. There are three scales involved in the analysis,
namely: the macroscopic or bulk mass, the microscopic size of the
fluid particle and the mesoscopic scale  corresponding to the pores.
In the references above one find descriptions of the velocity $v$ as
an average over the mesoscopic cells of the fluid particle
velocities. Taking into account that each cell contains a solid part
where the particle velocity vanishes, it is then natural to get
the viscous forces associated to that average velocity, which is a scaled
approximation of the laplacian term appearing in the Navier-Stokes
equation.

In this paper we shall consider the case of an homogeneous and
isotropic porous material. Porosity is the fraction of the volume
occupied by pore or void space. But it is important to distinguish
between two kind of pore, one that form a continuous interconnected
phase within the medium and the other consisting on isolated pores,
because non interconnected pores can not contribute to fluid
transport. Permeability is the term used to describe the
conductivity of the porous media with respect to a newtonian fluid
and it will depend upon the properties of the medium and the fluid.
Darcy's law indicates such dependence allowing us to define the notion
of specific permeability $\kappa$ and its appropriated units. In the
case of anisotropic material $\kappa$ will be a symmetric and
positive definite tensor, and then the methods of our proof can be
 modified to get local existence, but for a non
homogeneous medium the properties of the tensor $\kappa(x)$ will
have to be conveniently specified in order to have an interesting
theory.

The Muskat and related problems \cite{S-T} have been recently studied \cite{Peter,SCH,DP,DP2,ADP}. The first natural question asks for the evolution (existence) of such system,
at least for a short time $t>0$, and the persistence of smoothness of the interface $S(t)$
 if we begin with a smooth enough surface at time $t=0$. One can deduce easily from this
 formulation that in the occurrence of such smooth evolution both pressures, modulo a constant, must coincide at the interface:
$$
p^1|_{S(t)} = p^2|_{S(t)}.
$$
Therefore we look at the case without surface tension
(see  article \cite{ES} where the regularizing effect of surface tension is considered). The normal component of the velocity fields must also agree at the free boundary:
$$
(v^1 - v^2)\cdot\nu^j = 0\quad\quad \text{at} \quad S(t)
$$
where $\nu^j$ is the inner unit normal to $S$ at the domain $D^j$ ($\nu^2 = -\nu^1$).
Furthermore the vorticity will be concentrated at the interface, having form
$$
\curl (v) = \omega (z)dS(z)
$$
where $\omega$ is tangent at $S$ at the point $z$ and $dS(z)$ is surface measure.

The main purpose of this paper is to extend to the 3-dimensional case the results obtained
in \cite{ADP} for the case of 2 dimensions, namely proving local-existence in the scale of Sobolev
spaces of the initial value problem if the Rayleigh-Taylor condition (R-T) is initially satisfied (see \cite{S-T}) where this issue is studied from a physical point of view). In our case that condition amounts to the positivity of the function
$$
\sigma = (\nabla p^2 - \nabla p^1)\cdot (\nu^2-\nu^1)
$$
at the interface $S$. Let us indicate that the R-T property also appears in other fluid interface problems such as water waves \cite{ADP2}.

Together with that hypothesis, one also
assume  that the initial surface $S$ is connected and simply
connected. In the presence of a global  parametrization $X:
\mathbb{R}^2\rightarrow S$, the preservation of that character will
be controlled by the gauge
$$
F(X)(\al,\beta)=\frac{|\alpha - \beta|}{|X(\alpha) - X(\beta)|},\quad \|F(X)\|_{L^\infty}=
\sup_{\alpha\neq\beta}\frac{|\alpha - \beta|}{|X(\alpha) - X(\beta)|} <\infty.
$$

Section 2 of this paper contains the deduction of the evolution equations for the interface $S$.
 In section 3 we prove the existence of global isothermal parametrization as a consequence of
  the Koebe-Poincare uniformization theorem of Riemann surfaces in the geometric scenarios
  considered in our work, namely: either double periodicity in the horizontal variables or
   asymptotic flatness. Let us add that given the non-local character of the operator involved,
    to obtain a global isothermal parametrization is an important step in the proof, whose main components are sketched in section 4.

In closing our system (section 2) we need to control the norm of the inverse operator
$(I + \lambda \mathcal{D})^{-1}$ where $\mathcal{D}$ is the double-layer potential and
 $|\lambda|\leq 1$. It is well known from Fredholm\'{}s theory that those operators are
 bounded on $L^2(S)$. However since the surface $S = S(t)$ is moving, a precise control
  of its norm is needed in order to proceed with our proof. That is the purpose of section
  5 where the estimates for the double-layer potential are revisited.

In sections 6 and 7 we develop the energy estimates needed to conclude local-existence. Let us mention that at a crucial point (more precisely just at that step where the positivity of $\sigma (\alpha,t)$ (R-T) plays its role), we use the pointwise estimate $\theta(x)\Lambda\theta(x) \geq \frac{1}{2} \Lambda\theta^2(x)$ of \cite{CC}, with $\Lambda=\sqrt{-\Delta}$.

  In the strategy of our proof it is crucial to analyze the evolution of both quantities $\sigma$ and $F$ (section 8) at the same time than the interface $X$ and vorticity $\omega$. There are several publications (see \cite{Ambrose} for example) where different authors have treated these problems assuming that the Rayleigh-Taylor condition is preserved during the evolution. Under such hypothesis the proof can be considerably simplified, specially if one also assume the appropriated bounds for the resolvent of the double layer potential respect to a moving domain, or the existence of global isothermal coordinates, etc... It is our purpose of going carefully over such items what is responsible for the more delicate and intricate parts of this paper.

\section{The contour equation}

We consider the following evolution problem for the active scalars
$\rho = \rho(x,t)$ and $\mu=\mu (x,t)$, $x\in\mathbb{R}^3$, and $t\geq 0$:
\begin{equation*}
\rho_t + v\cdot\nabla\rho = 0,
\end{equation*}
\begin{equation*}
\mu_t + v\cdot\nabla\mu = 0,
\end{equation*}
with a velocity $v = (v_1,v_2,v_3)$ satisfying the momentum equation
\begin{equation}\label{dlaw}
\mu v=-\nabla p-(0,0,\,\rho),
\end{equation}
and the incompressibility condition $\nabla\cdot v = 0$, where, without loss of generality,   we have prescribed the values $\kappa=\mathrm{g}=1$.

The vector $(\mu,\rho)$ is defined by

\begin{equation*}
(\mu,\rho)(x_1,x_2,x_3,t)=\left\{\begin{array}{cl}
                    (\mu^1,\rho^1),& x\in D^1(t)\\
                    (\mu^2,\rho^2),& x\in D^2(t)=\R^3\setminus D^1(t),
                 \end{array}\right.
\end{equation*} where $\mu^1\neq\mu^2$, and $\rho^1\neq\rho^2$. Darcy's law \eqref{dlaw} implies that  the fluid is irrotational
in the interior of each domain $D^j$ and because of the jump of densities and viscosities on the free
boundary, we may assume a velocity field such that
$$
\curl v=\omega(\al,t)\delta(x-X(\al,t)),
$$
where $\partial D^j(t)=\{X(\al,t)\in \R^3: \al=(\al_1,\al_2)\in \R^2\}$, i.e.
\begin{equation}\label{defav}
<\curl v,\varphi>=\int_{\R^2} \omega(\al,t)\cdot\varphi(X(\al,t))d\al,
\end{equation}
for any $\varphi:\R^3\to\R^3$ vector field in $C^{\infty}_c(\R^3)$.

The incompressibility hypothesis ( $<\grad\cdot v,\varphi>\equiv -< v,\grad\varphi>=0$ for
any $\varphi\in C^{\infty}_c(\R^3)$) yields
$$
v^1(X(\al,t),t)\cdot N(\al,t) =v^2(X(\al,t),t)\cdot N(\al,t),
$$
with $N(\al,t)=\dau X(\al,t)\wedge \dad X(\al,t)$,
and equation \eqref{defav} gives us the identity
$$
\omega(\al,t)=(v^2(X(\al,t),t)-v^1(X(\al,t),t))\wedge N(\al,t).
$$
Defining the potential $\phi$ by
$v(x,t)=\grad \phi(x,t)$ for $x\in \R^2\setminus \partial D^j(t)$, we get
$$
\Omega(\al,t)=\phi^2(X(\al,t),t)-\phi^1(X(\al,t),t),
$$
$$\dau\Omega(\al,t)=(v^2(X(\al,t),t)-v^1(X(\al,t),t))\cdot \dau X,$$
$$\dad\Omega(\al,t)=(v^2(X(\al,t),t)-v^1(X(\al,t),t))\cdot \dad X.$$

Then, one has the equality
\begin{align*}
\omega(\al,t)&=(v^2(X(\al,t),t)-v^1(X(\al,t),t))\wedge (\dau X(\al,t)\wedge \dad X(\al,t))
\end{align*}
and therefore
\begin{equation}\label{oOf}
\omega(\al,t)=\dad\Omega(\al,t) \dau X(\al,t)-\dau\Omega(\al,t) \dad X(\al,t),
\end{equation}
implying that $\grad\cdot \curl v=0$ in a weak sense.

Using the law of Biot-Savart we have for $x$ not lying in the free surface ($x\neq X(\al,t)$) the following expression for the velocity:

$$
v(x,t)=-\frac{1}{4\pi}\int_{\R^2}\frac{x-X(\beta,t)}{|x-X(\beta,t)|^3}\wedge \omega(\beta) d\beta.
$$

It follows that
\begin{equation}\label{eqevcon}
X_t(\al)=BR(X,\omega)(\al,t)+C_1(\al)\dau X(\al)+C_2(\al)\dad X(\al),
\end{equation}
where $BR$ is the well-known Birkhoff-Rott integral
\begin{equation}\label{BR}
BR(X,\omega)(\al,t)=-\frac{1}{4\pi}PV\int_{\R^2}\frac{X(\al)-X(\beta)}{|X(\al)-X(\beta)|^3}\wedge \omega(\beta) d\beta.
\end{equation}
Next we will close the system using Darcy's law:

Since
$$
\grad \phi =v(x,t) -\Omega(\al,t) N(\al,t) \delta (x-X(\al,t))
$$
we have
$$
<\Delta\phi,\varphi>=-<\grad\phi,\grad\varphi>=\int_{\R^2} \Omega(\al,t)N(\al,t)\cdot \grad \varphi(X(\al,t))d\al,
$$
taking $\varphi(y)=-1/(4\pi |x-y|)$ one  obtain $\phi$ in terms of the double
layer potential:
$$
\phi(x)=-\frac{1}{4\pi}\int_{\R^2} \frac{x-X(\al)}{|x-X(\al)|^3}\cdot N(\al)\Omega(\al) d\al.
$$
Darcy's law yields
\begin{equation*}
\Delta p(x,t)=-\di (\mu(x,t)v(x,t))-\dxt \rho(x,t),
\end{equation*}
that is
\begin{equation*}
\Delta p(x,t)=P(\al,t)\delta (x-X(\al,t)),
\end{equation*}
where $P(\al,t)$ is given by
$$P(\al,t)=(\mu^2-\mu^1)v(X(\al,t),t)\cdot N(\al,t)+(\rho^2-\rho^1) N_3(\al,t),$$
implying the continuity of the pressure on the free boundary.

Next if $x\neq X(\al,t)$, i.e. $x$ is not placed at the interface, we can write  Darcy's law in the form
$$
\mu \phi(x,t)=-p(x,t)-\rho x_3
$$
and taking limits in both domains $D^j$ we get at $S$ the equality
$$(\mu^2\phi^2(X(\al,t),t)-\mu^1\phi^1(X(\al,t),t))=-(\rho^2-\rho^1)X_3(\al,t).$$
Then the formula for the double layer potential gives
$$
\frac{\mu^2+\mu^1}{2}\Omega(\al,t)-(\mu^2-\mu^1)
\frac{1}{4\pi}PV\int_{\R^2} \frac{X(\al)-X(\beta)}{|X(\al)-X(\beta)|^3}\cdot N(\beta)\Omega(\beta) d\beta=
-(\rho^2-\rho^1)X_3(\al,t)
$$
that is
\begin{equation}\label{eqOmega}
\Omega(\al,t)-A_\mu\mathcal{D}(\Omega)(\al,t)=-2A_\rho X_3(\al,t),
\end{equation}
where
\begin{equation}\label{dlp}
\mathcal{D}(\Omega)(\al)=
\frac{1}{2\pi}PV\int_{\R^2} \frac{X(\al)-X(\beta)}{|X(\al)-X(\beta)|^3}\cdot N(\beta)\Omega(\beta) d\beta,
\end{equation}
$$
A_\mu=\frac{\mu^2-\mu^1}{\mu^2+\mu^1}\qquad\mbox{and}\qquad A_\rho=\frac{\rho^2-\rho^1}{\mu^2+\mu^1}.
$$
 And the evolution equation are then given by \eqref{oOf}-\eqref{dlp}, where the functions $C_1$ and $C_2$ will be chosen in the next section.

Furthermore, taking limits we get from Darcy's law the following two formulas:
\begin{equation}\label{fduO}
\dau\Omega(\al,t)+2A_\mu BR(X,\omega)(\al,t)\cdot\dau X(\al,t)=
-2A_\rho\dau X_3(\al,t),
\end{equation}
\begin{equation}\label{fddO}
\dad\Omega(\al,t)+2A_\mu BR(X,\omega)(\al,t)\cdot\dad X(\al,t)=
-2A_\rho\dad X_3(\al,t).
\end{equation}

\section{Isothermal parameterization: choosing the tangential terms}

Although the normal component of the velocity vector field is the
relevant one in the evolution of the interface, it is however very
important to choose an adequate parameterization in order to uncover
and handle properly the cancelations  contained in the equations of
motion. Fortunately for our task we can rely upon the ideas of H.
 Lewy \cite{lewy}, and many other authors, who discovered the convenience of  using isothermal coordinates in different P.D.E.  namely for understanding how a minimal surface leaves an obstacle, but also in several fluid mechanical problems.

Let us recall that  an isothermal parameterization must satisfy:
\begin{equation*}\label{iso}
|X_{\alu}(\al,t)|^2=|X_{\ald}(\al,t)|^2, \qquad X_{\alu}(\al,t)\cdot X_{\ald}(\al,t)=0,
\end{equation*}
for $t\geq 0$.

Next we define
\begin{align}
\begin{split}\label{C1}
C_1(\al) &=\frac{1}{2\pi}\int_{\R^2}\frac{\al_1-\beta_1}{|\al-\beta|^2}\frac{BR_{\beta_2}
\cdot X_{\beta_2}- BR_{\beta_1}\cdot X_{\beta_1}}{|X_{\beta_2}|^2}d\beta\\
&\quad
-\frac{1}{2\pi}\int_{\R^2}\frac{\al_2-\beta_2}{|\al-\beta|^2}\frac{BR_{\beta_1}\cdot
X_{\beta_2}+ BR_{\beta_2}\cdot X_{\beta_1}}{|X_{\beta_1}|^2}d\beta,
\end{split}
\end{align}
and
\begin{align}
\begin{split}\label{C2}
C_2(\al) &=\frac{-1}{2\pi}\int_{\R^2}\frac{\al_2-\beta_2}{|\al-\beta|^2}\frac{BR_{\beta_2}
\cdot X_{\beta_2}- BR_{\beta_1}\cdot X_{\beta_1}}{|X_{\beta_2}|^2}d\beta\\
&\quad
-\frac{1}{2\pi}\int_{\R^2}\frac{\al_1-\beta_1}{|\al-\beta|^2}\frac{BR_{\beta_1}\cdot
X_{\beta_2}+ BR_{\beta_2}\cdot X_{\beta_1}}{|X_{\beta_1}|^2}d\beta.
\end{split}
\end{align}
That is $X_t=BR+C_1X_{\alu}+C_2X_{\ald}$ and
\begin{equation*}
X_{\alu t}=BR_{\alu}+C_1X_{\alu\alu}+C_2X_{\alu\ald}+C_{1\alu}X_{\alu}+C_{2\alu}X_{\ald},
\end{equation*}
\begin{equation*}
X_{\ald t}=BR_{\ald}+C_1X_{\alu\ald}+C_2X_{\ald\ald}+C_{1\ald}X_{\alu}+C_{2\ald}X_{\ald}.
\end{equation*}
Denoting $f=(|X_{\al_1}|^2-|X_{\al_2}|^2)/2$ and $g=X_{\al_1}\cdot X_{\al_2}$   we have
$$
f_t=(BR_{\al_1}\cdot X_{\al_1}-BR_{\al_2}\cdot X_{\al_2})+C_1
f_{\al_1}+C_2f_{\al_2}+(C_{2\al_1}-C_{1\al_2})g+ 2C_{1\al_1}f+
(C_{1\al_1}-C_{2\al_2})|X_{\al_2}|^2.$$

The expressions for $C_1$
and $C_2$ yield the vanishing of the sum of the first
and the last terms in the identity above. Therefore we get
\begin{equation}\label{elf}
 f_t=C_1 f_{\al_1}+C_2f_{\al_2}+(C_{2\al_1}-C_{1\al_2})g+2C_{1\al_1}f.
\end{equation}

Similarly we have
$$
g_t=(BR_{\al_2}\cdot X_{\al_1}+BR_{\al_1}\cdot X_{\al_2})+C_1
g_{\al_1}+C_2g_{\al_2}+(C_{1\al_1}+
C_{2\al_2})g-2C_{2\al_1}f+(C_{1\al_2}+C_{2\al_1})|X_{\al_1}|^2,
$$
and
\begin{equation}\label{elg}
g_t=C_1 g_{\al_1}+C_2g_{\al_2}+(C_{1\al_1}+C_{2\al_2})g-2C_{2\al_1}f.
\end{equation}

The linear character  of equations \eqref{elf} and \eqref{elg} allows us to conclude that if there is a solution
of the system $X_t=BR+C_1X_{\alu}+C_2X_{\ald}$ and we start with isothermal coordinates at time $t=0$,
then they will continue to be
isothermal so long as the evolution equations provide us with a smooth enough interface.


The fact that one can always prescribe such coordinates at time $t=0$ follows from the following argument:
In the double periodic setting we have a $C^2$ simply connected surface, homeomorphic to the euclidean
plane $\mathbb{R}^2$, which, by the Riemann-Koebe-Poincare uniformization theorem, is conformally
 equivalent to either the Riemann sphere, the plane or the unit disc. The sphere is easily eliminated by
  compacity, but we can also rule out the unit disc because the double periodicity assumption in the
  horizontal variables imply the existence of an abelian discrete subgroup of rank two in the group of
  conformal transformations, and that event  can not happen in the case of the unit disc.

  Therefore we have an orientation preserving  conformal (isothermal) equivalence
$$
\phi: \mathbb{R}^2 \longrightarrow  S.
$$
Since $S$ is invariant under translations $\tau_{\nu}(x)= x + 2 \pi \nu$, $\nu \in \mathbb{Z}^2 \times \{0\}$
 it follows that $f_{\nu}(z) = \phi^{-1}\circ\tau_{\nu}\circ\phi(z)$ must be a diffeoholomorphism of
 $\mathbb{C}=\mathbb{R}^2$ and, therefore, it has to be of the form
$$
f_{\nu}(z) = a_{\nu} z + b_{\nu}, \quad\text{for certain}\quad a_{\nu}, b_{\nu}\in \mathbb{C}.
$$
Clearly the family $f_{\nu}$ is generated by $f_1 =f_{(1,0,0)}$, $f_2=f_{(0,1,0)}$. Let
$$
f_1(z)= a_1 z + b_1, \quad\quad f_2(z)= a_2 z + b_2
$$
we claim that $a_1=a_2=1$.
Suppose that $|a_1|<1$ then we get $f_1^n(z)=a_1^n z + b_1(1 + a_1 + ..+a_1^{n-1})$ a sequence
converging to $\frac{b_1}{1-a_1}$ contradicting the discrete character of the group action.
On the other hand, if $|a_1|>1$ then since
$$
f_1^{-1}(z) = f_{(-1,0,0)}(z) = \frac{z}{a_1} - \frac{b_1}{a_1}
$$
we get a contradiction with the sequence $f_1^{-n}(z)$. Therefore we must have $a_1= e^{2\pi i\theta}$
for some $0\leq\theta<1$. Assume that $0<\theta<1$ then
$$
f_1^{(n)}(z) = e^{2\pi i n \theta}z+ b_1(1 + e^{2\pi i\theta}+.....+ e^{2\pi i(n-1)\theta})=
e^{2\pi in\theta}z + b_1\frac{1 - e^{2\pi in\theta}}{1- e^{2\pi i\theta}}
$$
That is the sequence $f^n(z)$ is bounded, $|f^n(z)|\leq |z| + \frac{|b_1|}{\sin\pi\theta}$, and
therefore it contains a converging subsequence contradicting again the discrete character of
the action. That is, we must have $f_1(z) = z + b_1$ and,  similarly, $f_2(z)= z + b_2$, allowing
us to conclude easily the double periodicity of the isothermal parameterization $\phi$.

In the asymptotically flat case we start with an orientable simply connected surface $S$ so
that outside a ball $B$ in $\mathbb{R}^3$ it becomes the graph of a $C^2$-function $x_3 = \varphi(x_1,x_2)$
 satisfying that $|D^{\alpha}\varphi(x)| = o(|x|^{-N})$, for every N and $|\alpha|\leq 2$, in particular, the normal
 vector $\frac{(-\nabla\varphi, 1)}{\sqrt{1 + |\nabla\varphi|^2}}=\nu(x)$ is pointing out
 vertically $\frac{1}{\sqrt{1 + |\nabla\varphi|^2}}\gg \frac{1}{2}$ for $|x|$ big enough.

It is then well known that one can find isothermal coordinates whose first fundamental form
  $ \lambda(\alpha,\beta)(d\alpha^2 + d\beta^2)$ converge asymptotically to the identity.

Again by the uniformation theorem $S$ must be conformally equivalent to either $\mathbb{C}$
 or the unit disc. But since outside $B$ the surface $S$ is conformally equivalent to
  $\mathbb{C}-B\bigcap\{x_3=0\}$ it cannot be also conformally equivalent to $D-K$, for
  any regular compact set $K$ contained in the unit disc $D$, because the harmonic measure
   of the ideal boundary  is 1 in the case of $D$ and 0 for $\mathbb{R}^2$.


\section{Outline of the proof.}


The proof of local existence requires the following:

1) A connected and simply connected surface $S=S(t)$ parameterized
by isothermal coordinates
$$
X:\mathbb{R}^2\longrightarrow\mathbb{R}^3, \quad X=X(\alpha,t)
$$
with normal vector $N(\alpha,t)= X_{\alpha_1} \wedge X_{\alpha_2}$ and gauge
$$
F(X)(\al,\beta)=|\beta|/|X(\al)-X(\al-\beta)|,
$$
such  $\|F(X)\|_{L^{\infty}}<\infty$ and $\||N|^{-1}\|_{L^{\infty}}<\infty$.

2) The positivity of
\begin{align}
\begin{split}\label{R-T}
\sigma(\al,t)&=-(\grad p^2(X(\al,t),t)-\grad p^1(X(\al,t),t))\cdot N(\al,t)\\
    &=(\mu^2-\mu^1)BR(X,\omega)(\al,t)\cdot N(\al,t)+(\rho^2-\rho^1)N_3(\al,t),
\end{split}
\end{align}
where the last equality is a consequence of  Darcy's law after taking limits in both domains $D^j$. This is the Rayleigh-Taylor condition to be imposed at time $t=0$, being a part of the problem to prove that it remains true as time pass.

3) The estimates on the norm of $(I - \lambda\mathcal{D})^{-1}$, $|\lambda|<1$, $\mathcal{D}=$ double layer potential (section 5) allows us to obtain the inequalities:
\begin{equation*}
\|\Omega\|_{H^{k+1}}\leq  P(\|X\|^2_{k+1}+\|F(X)\|^2_{L^\infty}+\||N|^{-1}\|_{L^\infty}),
\end{equation*}
\begin{equation*}
\|\omega\|_{H^k}\leq  P(\|X\|^2_{k+1}+\|F(X)\|^2_{L^\infty}+\||N|^{-1}\|_{L^\infty}),
\end{equation*}
for $k\geq 3$, where $P$ is a polynomial function and the norm $\|\cdot\|_k$  is given by
\begin{equation*}
\|X\|_k=\|X_1-\al_1\|_{L^3}+\|X_2-\al_2\|_{L^3}+\|X_3\|_{L^2}+\|\grad(X-(\al,0))\|^2_{H^{k-1}},
\end{equation*}
as in \eqref{nl3l2s} below, and $\|\cdot\|_{H^j}$ denotes the norm in Sobolev's space $H^j$.

4) A control of the Birkhoff-Rott integral $BR(X,\omega)$:
\begin{eqnarray*}
\|BR(X,\omega)\|_{H^k}\leq C P(\|X\|^2_{k+1}+\|F(X)\|^2_{L^\infty}+\||N|^{-1}\|_{L^\infty}).
\end{eqnarray*}
for $k\geq 3$.

5) Energy estimates: The  properties of  isothermal parameterizations help us to reorganize the terms in such a way that
\begin{align*}
\frac{d}{dt}\|X\|^2_{k}(t)&\leq  P(\|X\|^2_{k}(t)+\|F(X)\|^2_{L^\infty}(t)+\||N|^{-1}\|_{L^\infty}(t))\\
&\quad\D-\sum_{i=1,2}\frac{2^{3/2}}{(\mu_1\!+\!\mu_2)}\,\int_{\R^2} \frac{\sigma(\al,t)}{|\grad X(\al,t)|^{3}}
 \dai^{k} X(\al,t)\cdot \la(\dai^{k} X)(\al,t) d\al,
\end{align*}
where $k\geq 4$, $|\grad X(\al)|^{3}=(|\dau X(\al)|^{2}+|\dad X(\al)|^{2})^{3/2}$ and
 $\Lambda=(-\Delta)^{1/2}=(R_1(\dau)+R_2(\dad)).$ Then the pointwise inequality
 $$\theta\Lambda(\theta)-\frac12\Lambda(\theta^2)\geq 0,$$
 together with the condition $\sigma>0$ allows us to get rid of the dangerous terms in the inequality above (i.e. those involving
 (k+1)-derivatives of X) to obtain the estimate
 \begin{align*}
\frac{d}{dt}\|X\|^2_{k}(t)&\leq  P(\|X\|^2_{k}(t)+\|F(X)\|^2_{L^\infty}(t)+\||N|^{-1}\|_{L^\infty}(t)).
\end{align*}

 6) Finally we need to control the evolution of $\|F(X)\|_{L^\infty}(t)$ and  $\inf(t)=\D\inf_{\al\in \R^2} \sigma(\al,t)$  which is obtained via the following estimates
 \begin{equation*}
\D\dt\|F(X)\|^2_{L^\infty}(t)\leq  P(\|X\|^2_{4}(t)+\|F(X)\|^2_{L^\infty}(t)+\||N|^{-1}\|_{L^\infty}(t))
\end{equation*}
$$\frac{d}{dt}(\frac{1}{\inf(t)})\leq  \frac{1}{(\inf(t))^2}P(\|X\|^2_{4}(t)+\|F(X)\|^2_{L^\infty}(t)+\||N|^{-1}\|_{L^\infty}(t)).$$

7) All those facts together yield the inequality
$$\frac{d}{dt}E(t)\leq C P(E(t)),$$
for the energy:
$$
E(t)=\|X\|^2_{k}(t)+\|F(X)\|^2_{L^\infty}(t)+\||N|^{-1}\|_{L^\infty}(t)+(\inf(t))^{-1}
$$
where $k\geq 4$, $C$ is an universal constant and $P$ has polynomial growth
(depending upon $k$).

  At this point  it is not difficult to prove existence of a solution, locally in time, so long as the initial data $X(0)$ is in the appropriate Sobolev
space of order $k$, $k\geq 4$, and the Rayleigh-Taylor and not-selfintersecting conditions ($\sigma_0 > c >0, \|F(X(0))\|_{L^\infty} < \infty$) are satisfied.
Finally let us point out that since our existence proof is based upon energy inequalities an extra argument is needed to prove uniqueness. Nevertheless that task is much easier than proving existence ( the interested reader may consult the forthcoming paper \cite{ADP3} where the details of the proof have been written for several important cases, namely, Muskat, Water waves and SQG patches).

Let us remark that, at the end, we have to work with a coupled system involving the evolution of the surface $X$, the "vorticity density"
$\omega$, the Rayleigh-Taylor condition $\sigma$, the non-selfintersecting character of $S$ quantified by the gauge $F(X)$ and
the tangential parts $C_1 X_{\alpha_1} + C_2 X_{\alpha_2}$ of the velocity field.

This paper is a continuation of \cite{ADP} where the two-dimensional case was considered. Many of the needed estimates can be obtained following exactly the same methods that where used in \cite{ADP} for the lower dimensional case. Therefore, in order to simplify our presentation, we shall avoid here many details which were carefully proven in that quoted paper. This is specially the case of section 6 (control of the Birkhoff-Rott integral), section 8 (energy estimates) and also for the approximation schemes which are identical to those developed in \cite{ADP}. Therefore in the following we shall focus our attention on  the   more innovative parts of the proof, namely the evolution of the Rayleigh-Taylor condition, the non-selfintersecting property of the free boundary  and the needed estimates for double layer potentials.


\section{Inverting the operator: The single and double layer potentials revisited}

Along this proof we need to consider the properties of  single and double layer potentials,
which are well-known characters  in finding
solutions to the Dirichlet and Neumann problems in domains $D$ of $\mathbb{R}^n$.

For our purposes those domains will be of three different types,
namely:  bounded, periodic in the "horizontal" variables or
asymptotically flat. We shall also assume that their boundaries are
smooth enough (says $C^{2}$) and do not present self-intersections.
Therefore one has tangent balls at every point of the boundary, one
completely contained in $D$ and the other in $D^c$. We shall denote
by $\nu(x)$ the unit inner normal at the point $x\in\partial D$,
then under our hypothesis we have that, for $r>0$ small enough, the
parallel surfaces $\partial D_r = \{x + r\nu(x)| x\in\partial D\}$
are also $C^{2}$ surfaces with curvatures controlled by those of
$\partial D$. Furthermore the vector field $\nu$ can be extended
smoothly up to a collar neighborhood of $\partial D$ allowing us to
write the following formula:
\begin{eqnarray*}
\Delta u (x)= \frac{\partial^2 u}{\partial\nu^2}(x) - h(x)\frac{\partial u}{\partial\nu}(x) + \Delta_s u(x)
\end{eqnarray*}
where $\Delta$ denotes the ordinary laplacian in $\mathbb{R}^n$, $\Delta_s$ is the Laplace-Beltrami
operator in $\partial D$, $h(x)$ is the mean curvature of $\partial D$ at the point x and $u$ is any
$C^2$-function defined in a neighborhood of $\partial D$.

For convenience we will use the notation $D_1=D$, $D_2=D^c$,
$S=\partial D_j$ and $\nu_j(x) (j=1,2)$ the inner normal at $x\in S$
pointing inside $D_j$. Let $dS$ be the surface measure in $S$
induced by Lebesgue measure in ambience space, then given integrable
functions $\varphi, \psi$ on $S$ we have the integrals
$$
V(x) = c_n\int_S\psi(y)\frac{1}{||x-y||^{n-2}}dS(y)
$$

$$
W(x)=
c_n\int_S\varphi(y)\frac{\partial}{\partial\nu_x}(\frac{1}{||x-y||^{n-2}})dS(y)
$$
representing the single (respect. double) layer potential of $\psi$ (respect. $\varphi$), where $c_n$ is
 the normalizing constant so that $\frac{c_n}{||x||^{n-2}}$ becomes a fundamental solution of $\Delta$
 in $\mathbb{R}^n$, $n\geq 3$.

For $x\in S$ let us denote $W_1(x), V_1(x)$ (resp. $W_2(x), V_2(x)$) the corresponding limits
of the potentials inside $D_1$ (resp. $D_2$), we have:
$$
W_1(x)= \frac{1}{2}(\varphi(x) - \int_S \varphi(y) K(x,y)d\sigma(y))
= \frac{1}{2}(\varphi(x) - \mathcal{D}\varphi(x))
$$
$$
W_2(x)= \frac{1}{2}(\varphi(x) + \int_S \varphi(y) K(x,y)d\sigma(y))
= \frac{1}{2}(\varphi(x) + \mathcal{D}\varphi(x))
$$

$$
\frac{\partial V}{\partial\nu_1}(x) =-\frac{1}{2}(\psi(x) + \int_S
\psi(y) K(y,x)d\sigma(y)) =
- \frac{1}{2}(\psi(x) + \mathcal{D}^{\ast}\psi(x))\\
$$
$$
\frac{\partial V}{\partial\nu_2}(x) =-\frac{1}{2}(\psi(x) - \int_S
\psi(y) K(y,x)d\sigma(y)) = - \frac{1}{2}(\psi(x) -
\mathcal{D}^{\ast}\psi(x))
$$
where $$K(x,y)= 2c_n\frac{\partial}{\partial\nu_y}(\frac{1}{||x-y||^{n-2}}) =
\widetilde{c}_n\frac{\langle x-y,\nu(y)\rangle}{|x-y|^n}.$$

It is well known that in those scenarios considered above  the
boundary operators $\mathcal{D}$ (and $\mathcal{D}^{\ast}$)are
 smoothing of order $-1$ and therefore compact. Furthermore all their eigenvalues are real numbers
 having absolute value strictly less than $1$. Therefore, by the standard Fredholm theory, the
  operators $I-\lambda \mathcal{D}$, $I-\lambda \mathcal{D}^{\ast}$ are invertible when $|\lambda|\leq 1$. However,
  in our case the domains are moving and the evolution of their common boundary $S$ involves such
  inverse operators, making it necessary to estimate their norms in terms of the geometry and smoothness of $S$.

Although there is a vast literature about  single and double layer potentials, we have not been able to point out
 a precise statement giving the information needed for our results. Therefore in this section we provide arguments to
 prove that the norms of such inverse operators growth at most polynomially  $P(|||S|||)$, where $|||S|||$ is just
  $||S||_{C^{2}}$ plus a term of chord-arc type
  controlling the non-self-intersecting character of the boundary, namely we add the term $r(S)^{-1}$, where
   $r(S)$ is the $\sup$ over all the positive $r$ so that $S$ admits tangent  balls of radius $r$ in both domains
$D_j$: $$|||S||| = ||S||_{C^{2}} + (r(S))^{-1}.$$
   We shall write $P(|||S|||)$  to denote $\leq C(|||S|||^p)$
   for certain positive constants $C, p$ which are independent of the characters whose evolution is
   being controlled, but the size of both constants may change along the proof and we  shall make
   no effort to obtain their best values.

 We will consider the case of bounded domains in $\mathbb{R}^n$, $n\geq 3$, because the needed modifications when $n=2$, namely taking $\log |x|$ as fundamental solution for the laplacian, as well as the changes for the periodic or
 asymptotically flat domains, are left to the reader.

Let $\mathcal{D}$ and $\mathcal{D}^{\ast}$ be the potential defined
above with kernel
$$
K(x,y)= c_n\frac{\partial}{\partial\nu(y)}\frac{1}{||x - y||^{n-2}}= c_n \frac{\langle x-y,\nu(y)\rangle}{|x-y|^n}
$$
and $K(y,x)$ respectively. In the study of the inverse operators
$(I-\lambda \mathcal{D})^{-1}$, $|\lambda|\leq 1$
 it is convenient to consider first the particular values $\lambda = \pm 1$.

\begin{prop}\label{prop1} The following estimate holds
$$
||(I\pm \mathcal{D})^{-1}||_{L^2(S)}= P(|||S|||).
$$
\end{prop}

Since the boundedness of $(I\pm \mathcal{D})^{-1}$ in $L^2(S)$ is well known from
the general theory, we can simplify the proof considering only functions $f\in L^2(S)$
whose support lies inside a region of $S$ where the normal $\nu(x)$ is close enough
to a fixed direction. Then for a general $f$ an appropriate partition of unity would
 allow us to add the local estimates, so long as the number of pieces is controlled by $|||S|||$. We shall
 use the following observation of Rellich (lemma \ref{rellich}) whose proof is immediate

\begin{lemma} \label{rellich}. Let $u$ be a harmonic function and $h$ a smooth vector field in the domain $D$,
then we have: \begin{align*}i)&\quad \di(|\nabla u|^2 h) = 2 \di((\nabla u\cdot h)\nabla u) + O(|\nabla u|^2 |\nabla h|),\\
ii)&\quad\int_{\partial D}\langle \nu, h\rangle |\nabla u|^2 d\sigma = 2\int_{\partial D}\frac{\partial u}{\partial\nu} (\nabla u\cdot h)d\sigma +
O(\int_{ D}|\nabla u|^2 |\nabla h|).
\end{align*}
\end{lemma}

Given a function $f\in C^1(S)$ we may define $\nabla_{\tau} f$ choosing at each point $x\in S$ an
orthonormal basis $\{e_1,....,e_{n-1}\}$ of the tangent space $T_x(S)$ (we can consider also $\nabla_{\tau}f$ to
 be the gradient naturally associated to the induced Riemannian metric by the ambience space). In both ways, although
  different, we have that $|\nabla_{\tau}f|\equiv \Lambda_{\tau}f$ is an elliptic pseudo-differential
  operator of order $1$ in $S$.
Solving the Dirichlet problem  $\Delta u =0$ in $D$, $u|_S=f$ we obtain the operator
 $D_{\nu}\equiv \frac{\partial u}{\partial\nu}|_S$ which is also a pseudo-differential
 operator of order $1$ in $S$.

\begin{lemma}\label{number2}Let $f\in L^2(S)$ having support on the region
$\frac{1}{2}\leq\langle\nu(x),\eta\rangle\leq 1$ (for a fixed unit vector $\eta$), then we have:
$$
\int_S|D_{\nu}f|^2d\sigma\simeq \int_S|\nabla_{\tau}f|^2d\sigma
$$
where the constants involved in the stated equivalence $\simeq$ are  $P(|||S|||)$.
\end{lemma}
Proof: Let $u$ be harmonic in $D$ so that $u|_S=f$. Under our hypothesis about $f$ and
since $|\nabla u|^2 = |D_{\nu}u|^2 + |\nabla_{\tau}u|^2$ and $\nabla_{\tau}u$ is a
local operator $(\supp_S(\nabla_{\tau}f)\subset \supp(f))$, lemma \ref{rellich}  yields:
\begin{eqnarray*}
\frac{1}{2}\int_S|\nabla_{\tau}f|^2d\sigma \leq \int_S\langle\nu(x),\eta\rangle|\nabla_{\tau}u|^2d\sigma
\leq 3 \int_S|D_{\nu}u|^2d\sigma + 2 \int_S|\nabla_{\tau}u||D_{\nu}u|d\sigma
\end{eqnarray*}
from which we easily obtain
\begin{eqnarray*}
\int_S|\nabla_{\tau}f|^2d\sigma \leq P(|||S|||)\int_S|D_{\nu}f|^2d\sigma.
\end{eqnarray*}

To get the opposite inequality we proceed as before, but since $D_{\nu}f$ is not local, an extra argument is needed
to control the contribution of the region outside $\supp(f)$. Let us introduce surface discs $B_r(x)=\{y \in S|\|x-y\|\leq r\}$,
$x\in S$, $0\leq r\leq |||S|||^{-1}$ and domains $\Delta_r(x)=\{y + \rho\nu(x)| y\in B_r(x), \rho\leq r\}$. Given $R=\frac12 |||S|||^{-1}$
there exists a fixed unit vector $\eta$ so that $\frac{1}{2}\leq\langle\nu(y),\eta\rangle\leq 1$, for every $y\in B_R(x)$ and also a
smooth vector field $h$ such that $h|_{\Delta_R(x)}\equiv\eta$, $\supp(h)\subset\Delta_{2R}(x)$ and $\frac12|h(x)|\leq\langle h(x),\nu(x)\rangle$,
 $ ||\nabla h||^2\leq P(|||S|||)||h||$.

 In order to obtain the estimate
$$
\int_S|D_{\nu}f|^2d\sigma\leq P(|||S|||)\int_S|\nabla_{\tau}f|^2d\sigma
$$
we may assume, without loss of generality, that $\supp(f)\subset B_R(x)$, for some $x\in S$, and then prove that
 $$
\int_{B_R(y_0)}|D_{\nu}f|^2d\sigma\leq P(|||S|||)\int_S|\nabla_{\tau}f|^2d\sigma,
$$
uniformly on $y_0\in S$.

With the vector field $h$ defined above in $\Delta_{2R}(y)$ let us apply Rellich's estimate to get
$$
\int_S|D_{\nu}f|^2 \langle h, \nu(x)\rangle d\sigma(x) =  \int_S \langle \nu, h\rangle |\nabla_{\tau}f|^2d\sigma
-2 \int_S D_{\nu}f \nabla_{\tau}f \cdot  hd\sigma + O(\int_{ D}|\nabla u|^2 |\nabla h|)
$$
where $u$ satisfies $\Delta u = 0$ in $D$, $u|_S=f$. We get easily
$$
\int_{B_R(y_0)}|D_{\nu}f|^2 \langle h, \nu(x)\rangle d\sigma(x) = O( \int_S  |\nabla_{\tau}f|^2d\sigma
+ \int_{ D}|\nabla u|^2 |\nabla h|dx).
$$
Then the proof  will be finished if we can show that
$$
 \int_{ D}|\nabla u|^2 |\nabla h|dx \leq P(|||S|||)\int_S|\nabla_{\tau}f|^2d\sigma.
 $$
 To see it let us consider the parallel surfaces $S_r=\{x + r\nu(x)| x\in S\}$ ($0\leq r\leq |||S|||$) and observe that
 $$
 \int_{S_r}u^2d\sigma_r\simeq \int_S u^2(x + r\nu(x))d\sigma
 $$
 and
 \begin{eqnarray*}
 \int_S [u^2(x + r\nu(x)) - u^2(x)]d\sigma(x) = \int_S\int_0^r\nabla u^2(x + t\nu(x))\cdot \nu(x)dtd\sigma\\
 =2\int_{L_r}u(y)\nabla u(y)\cdot\nu(y) \leq 2(\int_{L_r}u^2(y))^{\frac{1}{2}}(\int_{L_r}|\nabla u|^2(y))^{\frac{1}{2}}
\end{eqnarray*}
where $L_r= \{x + \rho\nu(x)| x\in S, 0\leq \rho\leq r\}$.

Taking $F(x + r\nu(x))= f(x)\mathcal{X}(x)$ ($\mathcal{X}$= smooth cut-off) as a comparation function, Dirichlet's principle and
Poincare's inequality give us the estimate
$$
\int_D|\nabla u|^2\leq\int_D|\nabla F|^2\leq C(\int_S|\nabla_{\tau}f|^2 + \int_S|f|^2)= O(\int_S|\nabla_{\tau}f|^2d\sigma).
$$

Therefore
$$
 \int_{S_r}u^2d\sigma_r\simeq \int_S u^2(x + r\nu(x))d\sigma\leq \int_Sf^2(x)d\sigma + (\int_{L_r}u^2(y))^{\frac{1}{2}}(\int_{S}|\nabla_{\tau} f|^2)^{\frac{1}{2}}.
 $$
An integration in $r$, $0\leq r \leq R=|||S|||^{-1}$ yields
$$
\int_{L_r}u^2dx\leq R[ \int_Sf^2(x)d\sigma  + (\int_{L_r}u^2(y))^{\frac{1}{2}}(\int_{S}|\nabla_{\tau} f|^2)^{\frac{1}{2}}.
$$
That is
$$
\int_{L_r}u^2dx\leq CR\int_{S}|\nabla_{\tau} f|^2d\sigma.
$$
To conclude let us observe that
\begin{align*}
\int_D|\nabla u|^2|\nabla h| &= \frac12\int_D\Delta u^2|\nabla h| = \frac12 \int_D(\Delta u^2|\nabla h|-u^2\Delta(|\nabla h|)) + \frac12 \int_Du^2(|\nabla h|)\\
&= \frac12\int_S u\frac{\partial u}{\partial \nu}\cdot |\nabla h|d\sigma - \frac12 \int_S f^2\frac{(|\nabla h|)}{\partial\nu} d\sigma +
\frac12\int_D u^2\nabla|h|\\
&\leq (\int_{S} f^2d\sigma)^{\frac{1}{2}}(\int|\frac{\partial u}{\partial\nu}|^2|\nabla h|^2d\sigma)^{\frac{1}{2}} + C\int_{S} f^2d\sigma
+ C\int_{L_R}u^2.
\end{align*}

\underline{Proof of Proposition \ref{prop1}:} As before let $f\in C^1(S)$, $\supp(f)\subset\mathrm{U}_0$
 and let $u$ be
its single layer potential:
$$
u(x)= c_n\int_S\frac{f(y)}{||x-y||^{n-2}}dS(y)
$$
Then taking derivatives on each domain $D_j$ with respect to the normal direction and evaluating at $S$
 we get:
$$
\frac{\partial u}{\partial\nu_1} =- \frac{1}{2}(f(x) +
\mathcal{D}^{\ast}f(x)),
$$
$$
\frac{\partial v}{\partial\nu_2} =- \frac{1}{2}(f(x) -
\mathcal{D}^{\ast}f(x)).
$$

By lemma \ref{number2} we know that
\begin{eqnarray*}
\int_S|\frac{\partial v}{\partial\nu_1}|^2d\sigma\simeq\int_S|\nabla_{\tau}v|^2d\sigma\simeq\int_S|
\frac{\partial v}{\partial\nu_2}|^2d\sigma
\end{eqnarray*}
where the constants involved in the equivalences $\simeq$ are all
controlled by above (respect. below) by  $P(|||S|||)$ (respect.
$1/P(|||S|||)$).

Since $\frac{\partial v}{\partial\nu_1} + \frac{\partial v}{\partial\nu_2}= -f$ those estimates imply that
$$
\min(||f - \mathcal{D}^{\ast}f||_2, ||f +
\mathcal{D}^{\ast}f||_2)\geq \frac{1}{P(|||S|||)}
$$
i.e. $||(I\pm \mathcal{D})^{-1}||= P(|||S|||)$. Then using an appropriate partition of unity, that estimate
extends to a general
 $f\in L^2(S)$. q.e.d.

Next we shall consider Sobolev spaces $H^s(S)$, $0\leq s\leq 1$,
defined in the usual manner (i.e. throughout local coordinates
charts). We have also elliptic pseudo-differential operator
$\Lambda^s = (-\Delta)^{\frac{s}{2}}$ in such a way that
$$
||f||_{H^s(S)}\simeq ||f||_{L^2} + ||\Lambda^{s}f||_{L^2}.
$$

Then $H^{-s}(S)\equiv (H^{s}(S))^{\ast}$ allows us to consider the negative case by duality under the pairing
$$
\int_S\phi\psi d\sigma, \quad \phi\in H^{-s}, \quad \psi\in H^{s}
$$
and
$$
||\phi||_{H^{-s}}= \sup_{||\psi||_{H^s}=1}\int_S\phi\psi d\sigma
$$

Since both $\mathcal{D}$ and $\mathcal{D}^{\ast}$ are compact and
smoothing operators of degree $-1$, the commutators
 $[\Lambda^s,\mathcal{D}]$, $[\Lambda^s,\mathcal{D}^{\ast}]$
are then bounded in $L^2(S)$ $(0\leq s\leq 1)$ with norms controlled
by $|||S|||$, allowing us to extend proposition \ref{prop1} to the
chain of Sobolev's spaces:
\begin{cor} The norm of the operators $(I\pm \mathcal{D})^{-1}$, $(I\pm \mathcal{D}^{\ast})^{-1}$ in the space
$H^s(S)$, $-1\leq s \leq 1$, is bounded by $P(|||S|||)$.
\end{cor}


\subsection{Estimates for  $(I + \lambda \mathcal{D})^{-1}$,   $|\lambda|\leq 1$.}

With the same notation used before we have:
$$
\frac{1-\lambda}{2}\frac{\partial V}{\partial\nu_1} +
\frac{1+\lambda}{2}\frac{\partial V}{\partial\nu_2}= -
\frac{1}{2}(\phi(x) - \lambda \mathcal{D}^{\ast}\phi(x)),
$$
$$
\frac{1+\lambda}{2}\frac{\partial V}{\partial\nu_1} +
\frac{1-\lambda}{2}\frac{\partial V}{\partial\nu_2}= -
\frac{1}{2}(\phi(x) + \lambda \mathcal{D}^{\ast}\phi(x)),
$$
where
$$
V(x)= c_n\int_S\frac{\phi(y)}{||x-y||^{n-2}}dS(y).
$$

Then the identity $\phi -\lambda \mathcal{D}^{\ast}\phi=0$ yields
$$
0=(1-\lambda)\int_{\partial D_1}\!\!V\frac{\partial V}{\partial\nu_1}dS
+
 (1+\lambda)\int_{\partial D_2}\!\!V\frac{\partial V}{\partial\nu_2}dS
= (1-\lambda)\int_{ D_1}\!\!|\nabla V|^2 + (1+\lambda)\int_{ D_2}\!\!|\nabla V|^2
$$
which implies $\phi\equiv 0$. Similarly for $\phi +\lambda
\mathcal{D}^{\ast}\phi=0$, $-1\leq\lambda \leq 1$.

Remark: This observation can be improved applying the following fact
(whose proof we skip because it will not be used in our theorem):
$$
\int_{ D_1}|\nabla u|^2 \simeq \int_{ D_2}|\nabla u|^2
$$
where, again, the $\simeq$ is controlled by $P(|||S|||)$. In
particular it implies that the spectral radius of the operators
$\mathcal{D}$, $\mathcal{D}^{\ast}$ is less than $1-(P(|||S|||))^{-1}$.

\begin{thm}\label{thm2}
The operator norms $||(I + \lambda \mathcal{D})^{-1}||_{H^s(S)}$,
$||(I + \lambda \mathcal{D}^{\ast})^{-1}||_{H^s(S)}$, $|s|\leq 1$, $|\lambda|\leq 1$,
are $P(|||S|||)$ (growth at most polynomially with $|||S|||$).
\end{thm}
Proof: The identity $(I - \mathcal{D})^{-1}(I -\lambda \mathcal{D})
= I + (1 - \lambda)(I - \mathcal{D})^{-1} \mathcal{D}$ shows that
the conclusion of the theorem follows easily when $|1 - \lambda|\leq
\frac{1}{P(|||S|||)}$ and similarly when $|1 + \lambda|\leq
\frac{1}{P(|||S|||)}$.

Therefore, without loss of generality, we may assume that
$$
1-|\lambda|\geq \frac{1}{P(|||S|||)}.
$$

Assume now that $\phi\in H^{-\frac{1}{2}}(S)$ satisfies that $||\phi||_{H^{-\frac{1}{2}}}=1$ and
$$
||\phi - \lambda \mathcal{D}^{\ast}\phi||_{H^{-\frac{1}{2}}}\leq
\frac{1}{P(|||S|||)}.
$$
Then the single layer potential
$$
V(x)= c_n\int_S\frac{\phi(y)}{||x-y||^{n-2}}dS(y)
$$
satisfies the inequality
$$
|\int_S V(\phi - \lambda \mathcal{D}^{\ast}\phi)dS |\leq
\frac{1}{P(|||S|||)}
$$
On the other hand one have
$$
\int_S V(\phi - \lambda \mathcal{D}^{\ast}\phi)dS =
(1-\lambda)\int_{D_1}|\nabla V|^2 + (1+\lambda)\int_{D_2}|\nabla
V|^2
$$
implying the estimate
$$
\int_S V(\phi + \lambda \mathcal{D}^{\ast}\phi)dS =
(1+\lambda)\int_{D_1}|\nabla V|^2 + (1-\lambda)\int_{D_2}|\nabla
V|^2\leq \frac{1}{P(|||S|||)}.
$$
Then adding both inequalities together we would obtain
$$
\int_S V\phi d\sigma\leq \frac{1}{P(|||S|||)}
$$
which is impossible because of the following:
\begin{lemma} \label{number3}If V is the single layer potential of
$\phi$ then
$$
\int_S V(x)\phi(x) dS(x)=\int_S\int_S
\frac{\phi(x)\phi(y)}{||x-y||^{n-2}}
dS(x)dS(y)\geq\frac{1}{P(|||S|||)}||\phi||^2_{H^{-\frac{1}{2}}(S)}.
$$
\end{lemma}
Proof: Let us observe first that
$$
\int_S\int_S \frac{\phi(x)\phi(y)}{||x-y||^{n-2}} d\sigma(x)d\sigma(y)=
\int_{\mathbb{R}^n}\frac{1}{|\xi|^2}|\widehat{\phi d\sigma(\xi)}|^2d\xi\geq 0,
$$
where $\widehat{\phi dS}$ denotes the Fourier transform of the
measure $\phi dS$
 supported on $S$. This implies that
$$
\langle\phi,\psi\rangle = \int_S\int_S
\frac{\phi(x)\phi(y)}{||x-y||^{n-2}} dS(x)dS(y)
$$
is an inner product satisfying:
$$
|\langle\phi,\psi\rangle|\leq \langle\phi,\phi\rangle^{\frac{1}{2}}\langle\psi,\psi\rangle^{\frac{1}{2}}
$$
and we wish to show that
$$
\langle\phi,\phi\rangle\backsimeq ||\phi||_{H^{-\frac{1}{2}}(S)}^2.
$$
To see it let us observe first that given $\phi\in
H^{-\frac{1}{2}}(S)$ then its single layer potential $u|S$ belong to
the space $H^{\frac{1}{2}}(S)$ satisfying:
$$
||u||_{H^{\frac{1}{2}}(S)}\leq
P(|||S|||)||\phi||_{H^{-\frac{1}{2}}(S)},
$$
which  can be proved easily using local coordinates. As a
consequence we have
$$
\int_S\int_S \frac{\phi(x)\phi(y)}{||x-y||^{n-2}} dS(x)dS(y)\leq
P(|||S|||)||\phi||_{H^{-\frac{1}{2}}(S)}^2.
$$
In the opposite direction, since
$H^{-s}=(H^s)^{\ast}$ we have
$$
||\phi||_{H^{-s}}=\sup_{f\in H^s}\int_S\phi(x)f(x)d\sigma(x).
$$
Let us assume, for the moment,  that given $f\in H^s$ there exists
$g\in H^{s-1}$ so that
$$
f(x)= c_n\int_S\frac{g(y)}{||x-y||^{n-2}}dS(y)
$$
and $||f||_{H^s}\backsimeq ||g||_{H^{s-1}}$ (where we have used again the
symbol $\backsimeq$ to denote equivalence modulo a factor
$P(|||S|||)$).

Then
$$
||\phi||_{H^{-s}}\backsimeq    \sup_{||g||_{ H^{s-1}}=1}\langle\phi,g\rangle,
$$
and taking $s=\frac{1}{2}$, $s-1=-\frac{1}{2}$ we get
$$
||\phi||_{H^{-\frac{1}{2}}}\leq P(|||S|||)\langle\phi,
\phi\rangle^{\frac{1}{2}}\langle g, g\rangle^{\frac{1}{2}}\leq P(|||S|||)\langle\phi,\phi\rangle^{\frac{1}{2}}
||g||_{H^{-\frac{1}{2}}}\leq P(|||S|||) \langle\phi,\phi\rangle^{\frac{1}{2}}.
$$

To close our argument it remains to solve the equation
$$
f(x)= c_n\int_S\frac{g(y)}{||x-y||^{n-2}}dS(y)
$$
i.e. to prove that given given $f\in H^s$ there exits $g\in H^{s-1}$
satisfying the relation above.

To see that let us consider the solution of the Dirichlet problem:
\begin{equation*}
\begin{cases}
  \Delta u=0 \quad in \quad D_1\\
u|_S=f
\end{cases}
\end{equation*}
and the equation
$$
-2\frac{\partial u}{\partial\nu_1}= g - \mathcal{D}^{\ast}g
$$
i.e. $g= (I- \mathcal{D}^{\ast})^{-1}(-2\frac{\partial
u}{\partial\nu_1})$. Then we claim that such $g$ verifies the identity
$$
f(x)= c_n\int_S\frac{g(y)}{||x-y||^{n-2}}dS(y).
$$
This is because the function
$$
V(x)= c_n\int_S\frac{g(y)}{||x-y||^{n-2}}dS(y)
$$
is harmonic in $D_1$ and satisfies
$$
-2\frac{\partial V}{\partial\nu_1}= g -
\mathcal{D}^{\ast}g=-2\frac{\partial u}{\partial\nu_1},
$$
which implies that $V=u$ in $D_1$ and, therefore, taking limits up
to the boundary we obtain
$$
f(x)= c_n\int_S\frac{g(y)}{||x-y||^{n-2}}dS(y).
$$

To finish the proof of theorem \ref{thm2} let us consider for every
$\tau$, $0\leq\tau\leq 1$, the identity
$$
(I - \lambda \mathcal{D})^{-1}\Lambda^{\tau} = \Lambda^{\tau}(I -
\lambda \mathcal{D})^{-1} +
 (I - \lambda \mathcal{D})^{-1}C_{\tau}(I - \lambda \mathcal{D})^{-1}
$$
where the commutator $C_{\tau}=[\mathcal{D}\Lambda^{\tau} -
\Lambda^{\tau} \mathcal{D}]$ is a pseudodifferential operator of
order $\tau-2$ whose bounds are controlled by $|||S|||$. Then
\begin{align*}
||(I - \lambda \mathcal{D})^{-1}f||_{H^s}&\leq ||(I - \lambda
\mathcal{D})^{-1}f||_{H^{-\frac{1}{2}}} +
||\Lambda^{s + \frac{1}{2}}(I - \lambda \mathcal{D})^{-1}f||_{H^{-\frac{1}{2}}}\\
&\lesssim ||f||_{H^{-\frac{1}{2}}} + ||(I - \lambda \mathcal{D})^{-1}\Lambda^{s + \frac{1}{2}}f||_{H^{-\frac{1}{2}}}\\
&\lesssim ||f||_{L^2} + ||\Lambda^{s +
\frac{1}{2}}f||_{H^{-\frac{1}{2}}}\leq P(|||S|||)||f||_{H^{s}}
\end{align*}
q.e.d.

\begin{rem}
In the particular case of the sphere $S=S^{n-1}$ ($n\geq 2$) the estimate of lemma 5.6 becomes an identity:
$$
\int_{S^{n-1}}\int_{S^{n-1}}\frac{\phi(x)\phi(y)}{||x - y||^{n-2}}dS(x)dS(y) = c_n ||\phi||^2_{H^{-\frac{1}{2}}(S^{n-1})}
$$
for $n\geq 3$, and
$$
-\int_{S^{1}}\int_{S^{1}}\log||x - y||\phi(x)\phi(y)dS(x)dS(y) = c_2||\phi||^2_{H^{-\frac{1}{2}}(S^{1})}
$$
for $n=2$.

Proof: We present the details when $n\geq 3$. The case $n=2$ follows similarly.
Let $\phi(x)=\sum a_kY_k(x)$ where $Y_k$ is a spherical harmonic of degree $k$ normalized so that $||Y_k||_{L^2(S^{n-1})}=1$
then we have
$$
|a_0|^2 + \sum_{k\geq 1}\frac{|a_k|^2}{2k + n -2}=||\phi||^2_{H^{-\frac{1}{2}}(S)}<\infty.
$$
Claim: If $k\neq j$ then
$$
\int_{S^{n-1}}\int_{S^{n-1}}\frac{Y_k(x)Y_j(y)}{||x - y||^{n-2}}dS(x)dS(y) = 0
$$
Taking the Fourier transform and using Plancherel we get
$$
\int_{S^{n-1}}\int_{S^{n-1}}\frac{Y_k(x)Y_j(y)}{||x - y||^{n-2}}dS(x)dS(y) = \int_{\mathbb{R}^n}\frac{1}{|\xi|^2}\widehat{Y_kdS(\xi)}\overline{\widehat{Y_jdS(\xi)}}d\xi
$$
But it turns out that
$$
\widehat{Y_kdS(\xi)}= 2\pi i^{-k}|\xi|^{\frac{n-2}{2}}J_{\frac{n+2k-2}{2}}(|\xi|)Y_k(\frac{\xi}{|\xi|})
$$
where $J_{\nu}$ designs Bessel's function of order $\nu$, implying the claim.

Therefore our estimate diagonalizes:
 $$
 \int_{\mathbb{R}^n}\frac{1}{|\xi|^2}|\widehat{Y_kdS(\xi)}|^2d\xi= c\int_0^{\infty}\frac{1}{r}|J_{k+\frac{n-2}{2}}(r)|^2dr
$$
and the following well known identity for Bessel's functions
$$
\int_0^{\infty}\frac{J_{\mu}^2(r)}{r}dr = \frac{1}{2}\frac{1}{\mu}
$$
allows us to finish the proof.
\end{rem}


\subsection{Estimates for $\Omega$ and $\omega$.}

In the following we shall consider asymptotically flat  domains leaving to the reader the details of the periodic case. Since  we have  controlled  the norms of the operator relating $\Omega$ and $X$, we are in position to obtain the  following inequality:
\begin{equation}\label{eO}
\|\Omega\|_{H^{k}}\leq  P(\|X\|^2_{k}+\|F(X)\|^2_{L^\infty}+\||N|^{-1}\|_{L^\infty}),
\end{equation}
for $k\geq 4$, with $P$ a polynomial function. Then Sobolev's embedding implies
\begin{equation}\label{eo}
\|\omega\|_{H^k}\leq  P(\|X\|^2_{k+1}+\|F(X)\|^2_{L^\infty}+\||N|^{-1}\|_{L^\infty}),
\end{equation}
for $k\geq 3$. We will present the proof \eqref{eO} when $k=4$, because the cases  $k>4$ can be obtained with the same method.

 Theorem \ref{thm2} in \eqref{eqOmega} yields
$$
\|\Omega\|_{H^1}=\|(I-A_\mu\mathcal{D})^{-1}(-2A_\rho X_3)\|_{H^1}\leq C\|(I-A_\mu\mathcal{D})^{-1}\|_{H^1}\|X_3\|_{H^1}\leq P(|||S|||)\|X_3\|_{H^1}
$$
implying that
$$
\|\Omega\|_{H^1}\leq P(\|X\|^2_{4}+\|F(X)\|^2_{L^\infty}+\||N|^{-1}\|_{L^\infty}).
$$
Next we will show that
\begin{equation}\label{h2O}
\|\dau^2\Omega\|_{L^2}\leq P(\|X\|^2_{4}+\|F(X)\|^2_{L^\infty}+\||N|^{-1}\|_{L^\infty})\|\Omega\|_{H^1}
\end{equation}
 which together with the estimate for $\|\Omega\|_{H^1}$ above will allows us to control $\dau^2\Omega$ in terms of the free boundary.

 In order to do that we start with  formula \eqref{fduO} to get $\dau^2\Omega=I_1+I_2+I_3+I_4-2A_\rho \dau^2X_3$
where
$$
I_1=\frac{A_\mu}{2\pi}PV\int_{\R^2}
\frac{X(\al)-X(\al-\beta)}{|X(\al)-X(\al-\beta)|^3}\wedge\omega(\al-\beta)d\beta\cdot\dau^2 X(\al).
$$
$$
I_2=\frac{A_\mu}{2\pi}PV\int_{\R^2}
\frac{\dau X(\al)-\dau X(\al-\beta)}{|X(\al)-X(\al-\beta)|^3}\wedge\omega(\al-\beta)d\beta \cdot\dau X(\al),
$$
$$
I_3=-\frac{3A_\mu}{4\pi}PV\int_{\R^2}
A(\al,\beta)\frac{X(\al)- X(\al-\beta)}{|X(\al)-X(\al-\beta)|^5}\wedge\omega(\al-\beta)d\beta \cdot\dau X(\al),
$$
with $A(\al,\beta)=(X(\al)- X(\al-\beta))\cdot(\dau X(\al)-\dau X(\al-\beta))$, and
$$
I_4=\frac{A_\mu}{2\pi}PV\int_{\R^2}
\frac{X(\al)-X(\al-\beta)}{|X(\al)-X(\al-\beta)|^3}\wedge\dau\omega(\al-\beta)d\beta\cdot\dau X(\al).
$$
Our next objective is to introduce the operators  $\mathcal{T}_k$ \eqref{operadorgeneral} defined at the  appendix, in the analysis of those integrals $I_j$. Formula \eqref{oOf} gives us  $\omega=\dad(\Omega\dau X)-\dau(\Omega\dad X)$ and from standard Sobolev's estimates we get
$$
\|I_j\|_{L^2}\leq P(\|X\|^2_{4}+\|F(X)\|^2_{L^\infty}+\||N|^{-1}\|_{L^\infty})\|\Omega\|_{H^1}, \quad j=1,2,
$$
and similarly with $I_3$.

Regarding
$$I_4=\int_{|\beta|>1}d\beta+\int_{|\beta|<1}d\beta=J_1+J_2$$
we integrate by parts in $J_1$ to obtain
\begin{align*}
J_1&=\frac{A_\mu}{2\pi}\int_{|\beta|>1}\dbu\Big(\frac{X(\al)-X(\al-\beta)}{|X(\al)-X(\al-\beta)|^3}\Big)\wedge\omega(\al-\beta)d\beta\cdot\dau X(\al)\\
&\quad-\frac{A_\mu}{2\pi}\int_{|\beta|=1}\frac{X(\al)-X(\al-\beta)}{|X(\al)-X(\al-\beta)|^3}\wedge\omega(\al-\beta)dl(\beta)\cdot\dau X(\al).
\end{align*}
 From this last expression it is easy to deduce the inequality
\begin{align*}
J_1&\leq C\|F(X)\|^3_{L^\infty}\|X-(\al,0)\|^2_{C^1}\big(\int_{|\beta|>1}\frac{|\omega(\al-\beta)|}{|\beta|^3}d\beta+\int_{|\beta|=1}|\omega(\al-\beta)|dl(\beta)\big)
\end{align*}
 providing us with an appropriated control  (see appendix for more details).

 Next let us consider
 $J_2=K_1+K_2+K_3+K_4$ where
$$
K_1=\frac{A_\mu}{2\pi}PV\int_{|\beta|<1}
\frac{X(\al)-X(\al-\beta)}{|X(\al)-X(\al-\beta)|^3}\wedge\dad\Omega(\al-\beta)\dau^2 X(\al-\beta)d\beta\cdot\dau X(\al),
$$
$$
K_2=\frac{A_\mu}{2\pi}PV\int_{|\beta|<1}
\frac{X(\al)-X(\al-\beta)}{|X(\al)-X(\al-\beta)|^3}\wedge\dau\dad\Omega(\al-\beta)\dau X(\al-\beta)d\beta\cdot\dau X(\al),
$$$$
K_3=-\frac{A_\mu}{2\pi}PV\int_{|\beta|<1}
\frac{X(\al)-X(\al-\beta)}{|X(\al)-X(\al-\beta)|^3}\wedge\dau\Omega(\al-\beta)\dau\dad X(\al-\beta)d\beta\cdot\dau X(\al),
$$$$
K_4=-\frac{A_\mu}{2\pi}PV\int_{|\beta|<1}
\frac{X(\al)-X(\al-\beta)}{|X(\al)-X(\al-\beta)|^3}\wedge\dau^2\Omega(\al-\beta)\dad X(\al-\beta)d\beta\cdot\dau X(\al).
$$
Then the terms  $K_1$ and $K_3$ are handled with  the same approach used for $I_2$ ( i.e. \eqref{dT1} in the  appendix) and
we rewrite $K_2$ in the form
$$
K_2=\frac{A_\mu}{2\pi}\!\!\int_{|\beta|<1}
\frac{X(\al)\!-\!X(\al\!-\!\beta)}{|X(\al)\!-\!X(\al\!-\!\beta)|^3}\wedge\dau\dad\Omega(\al\!-\!\beta)(\dau X(\al\!-\!\beta)\!-\!\dau X(\al))d\beta\cdot\dau X(\al),
$$
to show that it can be estimated via an integration  by parts in the variable $\beta_1$ using the identity
$$\dau\dad\Omega(\al\!-\!\beta)=-\dbu(\dad\Omega(\al\!-\!\beta))$$ and the fact
 that the kernel in the integral $K_2$ has degree $-1$.

 It remains to deal with $K_4$: to do that  let us  consider $K_4=L_1+L_2$ where
$$
L_1=\frac{A_\mu}{2\pi}PV\!\!\int_{|\beta|<1}
\frac{X(\al)\!-\!X(\al\!-\!\beta)}{|X(\al)\!-\!X(\al\!-\!\beta)|^3}\wedge\dau^2\Omega(\al\!-\!\beta)(\dad X(\al)\!-\!\dad X(\al\!-\!\beta))d\beta\cdot\dau X(\al),
$$
and
$$
L_2=\frac{A_\mu}{2\pi}PV\int_{|\beta|<1}
\frac{X(\al)-X(\al-\beta)}{|X(\al)-X(\al-\beta)|^3}\dau^2\Omega(\al-\beta)d\beta\cdot N(\al).
$$
The term $L_1$ can be controlled like $K_2$, and  $L_2$ can be rewritten in the form:
$$
L_2=\frac{A_\mu}{2\pi}PV\int_{|\beta|<1}\Big(
\frac{X(\al)-X(\al-\beta)}{|X(\al)-X(\al-\beta)|^3}-\frac{\grad X(\al)\cdot\beta}{|\grad X(\al)\cdot\beta|^3}\Big)\dau^2\Omega(\al-\beta)d\beta\cdot N(\al),
$$
 showing that it can be estimated as we did with $\mathcal{T}_4$ \eqref{operadorregular}, that is we obtain \eqref{h2O}.
Similarly, equation \eqref{fddO} yields
\begin{equation*}
\|\dad^2\Omega\|_{L^2}\leq P(\|X\|^2_{4}+\|F(X)\|^2_{L^\infty}+\||N|^{-1}\|_{L^\infty})\|\Omega\|_{H^1},
\end{equation*}
and then the inequality $2\|\dau\dad\Omega\|_{L^2}\leq \|\dau^2\Omega\|_{L^2}+\|\dad^2\Omega\|_{L^2}$ gives us the desired control upon $\|\Omega\|_{H^2}$.

Next we will show that
\begin{equation}\label{h3O}
\|\dau^3\Omega\|_{L^2}\leq P(\|X\|^2_{4}+\|F(X)\|^2_{L^\infty}+\||N|^{-1}\|_{L^\infty})\|\Omega\|_{H^2}
\end{equation}
 allowing us to use the estimates for $\|\Omega\|_{H^2}$ above. In order to do that we start with formula \eqref{fduO} to get $\dau^3\Omega=\dau I_1+\dau I_2+\dau I_3+\dau I_4-2A_\rho \dau^3X_3$ where the most singular terms are given by
$$
J_3=\frac{A_\mu}{2\pi}PV\int_{\R^2}
\frac{X(\al)-X(\al-\beta)}{|X(\al)-X(\al-\beta)|^3}\wedge\omega(\al-\beta)d\beta\cdot\dau^3 X(\al),
$$
$$
J_4=\frac{A_\mu}{2\pi}PV\int_{\R^2}
\frac{\dau^2 X(\al)-\dau^2 X(\al-\beta)}{|X(\al)-X(\al-\beta)|^3}\wedge\omega(\al-\beta)d\beta \cdot\dau X(\al),
$$
$$
J_5=-\frac{3A_\mu}{4\pi}PV\int_{\R^2}
B(\al,\beta)\frac{X(\al)- X(\al-\beta)}{|X(\al)-X(\al-\beta)|^5}\wedge\omega(\al-\beta)d\beta \cdot\dau X(\al),
$$
with $B(\al,\beta)=(X(\al)- X(\al-\beta))\cdot(\dau^2 X(\al)-\dau^2 X(\al-\beta))$, and
$$
J_6=\frac{A_\mu}{2\pi}PV\int_{\R^2}
\frac{X(\al)-X(\al-\beta)}{|X(\al)-X(\al-\beta)|^3}\wedge\dau^2\omega(\al-\beta)d\beta\cdot\dau X(\al),
$$
and where the remainder terms can be estimated with the same method used before.

 Now we  write $$J_3=\frac{A_\mu}{2\pi}\mathcal{T}_1(\dad(\Omega\dau X)-\dau(\Omega\dad X))\cdot \dau^3 X$$
to obtain:
$$\|J_3\|_{L^2}\leq C \|\mathcal{T}_1(\dad(\Omega\dau X)-\dau(\Omega\dad X))\|_{L^4}\|\dau^3 X\|_{L^4}.$$
Next let us observe that in the proof of  estimate \eqref{estimaciongeneral} one can replace $L^2$ by $L^p$ for $1<p<\infty$ (see \cite{St3}). In particular we have
$$
\|J_3\|_{L^2}\leq P(\|X-(\al,0)\|_{C^{1,\delta}}+\|F(X)\|_{L^\infty}+\||N|^{-1}\|_{L^\infty})(\|\Omega\dau X\|_{L^4}+\|\Omega\dad X\|_{L^4}+\|\omega\|_{L^4})\|\dau^3 X\|_{L^4},$$
and then  Sobolev's embedding in dimension two: ($\|g\|_{L^4}\leq C\|g\|_{H^1}$)
yields the desired control. Regarding  $J_4$ we follow the approach taken before for $\mathcal{T}_3$ but using now the $L^4$ norm. That is we split
$$
J_4=\int_{|\beta|>1}d\beta+\int_{|\beta|<1}d\beta=K_5+K_6
$$
and since
$$
K_5\leq \|X-(\al,0)\|^2_{C^{2}}\|F(X)\|^3_{L^\infty}\int_{|\beta|>1}\frac{|\omega(\al-\beta)|}{|\beta|^3}d\beta
$$
that term can be estimated as above.

Next we introduce the splitting $K_6=L_3+L_4$ where
$$
L_3\!=\!\frac{A_\mu}{2\pi}\!\int_{|\beta|<1}\!\!\!\!\!\!\!\!\!(\dau^2 X(\al)\!-\!\dau^2 X(\al\!-\!\beta))[\frac{1}{|X(\al)\!-\!X(\al\!-\!\beta)|^3}-\frac{1}{|\grad  X(\al)\!\cdot\!\beta|^3}]
\wedge\omega(\al\!-\!\beta)d\beta\!\cdot\!\dau X(\al),
$$
$$
L_4=\frac{A_\mu}{2\pi}PV\int_{|\beta|<1}
\frac{\dau^2 X(\al)-\dau^2 X(\al-\beta)}{|\grad X(\al)\cdot\beta|^3}\wedge\omega(\al-\beta)d\beta \cdot\dau X(\al).
$$
We have
$$L_3\leq C\|X-(\al,0)\|^3_{C^{2,\delta}}(\|F(X)\|^4_{L^\infty}+\|X-(\al,0)\|^4_{C^1}\||N|^{-1}\|^4_{L^\infty})\int_{|\beta|<1}\frac{|\omega(\al-\beta)|}{|\beta|^{2-\delta}}d\beta $$
(see appendix for more details), giving us the appropriated estimate. Regarding $L_4$ we use identity \eqref{lo} which after a careful  integration by parts yields
\begin{align*}
L_4&=\frac{A_\mu}{2\pi}PV\int_{|\beta|<1}
\frac{\beta\cdot\grad_{\beta}\big((\dau^2 X(\al)-\dau^2 X(\al-\beta))\wedge\omega(\al-\beta)\cdot\dau X(\al)\big)}{|\grad X(\al)\cdot\beta|^3}d\beta\\
&\quad -\frac{A_\mu}{2\pi}\int_{|\beta|=1}
\frac{|\beta|(\dau^2 X(\al)-\dau^2 X(\al-\beta))\wedge\omega(\al-\beta)\cdot\dau X(\al)}{|\grad X(\al)\cdot\beta|^3}dl(\beta).
\end{align*}
helping us to prove the inequality
$$
\|L_4\|_{L^2}\leq P(\|X-(\al,0)\|_{C^2}+\|F(X)\|_{L^\infty}+\||N|^{-1}\|_{L^\infty})
(\|\dau^3 X\|_{L^4}\|\omega\|_{L^4}+\|\omega\|_{L^2}).$$
 Clearly $J_5$ can be approached with the same method used for $J_4$. Regarding the term $J_6$ we have to decompose further: first its most singular terms which are given by
$$
L_5=\frac{A_\mu}{2\pi}PV\int_{|\beta|<1}
\frac{X(\al)-X(\al-\beta)}{|X(\al)-X(\al-\beta)|^3}\wedge\dad\Omega(\al-\beta)\dau^3 X(\al-\beta)d\beta\cdot\dau X(\al),
$$
$$
L_6=\frac{A_\mu}{2\pi}PV\int_{|\beta|<1}
\frac{X(\al)-X(\al-\beta)}{|X(\al)-X(\al-\beta)|^3}\wedge\dau^2\dad\Omega(\al-\beta)\dau X(\al-\beta)d\beta\cdot\dau X(\al),
$$$$
L_7=-\frac{A_\mu}{2\pi}PV\int_{|\beta|<1}
\frac{X(\al)-X(\al-\beta)}{|X(\al)-X(\al-\beta)|^3}\wedge\dau\Omega(\al-\beta)\dau^2\dad X(\al-\beta)d\beta\cdot\dau X(\al),
$$$$
L_8=-\frac{A_\mu}{2\pi}PV\int_{|\beta|<1}
\frac{X(\al)-X(\al-\beta)}{|X(\al)-X(\al-\beta)|^3}\wedge\dau^3\Omega(\al-\beta)\dad X(\al-\beta)d\beta\cdot\dau X(\al).
$$
 Second let us observe  that the remainder is easy to handle; the terms $L_5$ and $L_7$ can be estimated as we did with $K_1$ and $K_3$ using the $L^4$ norm and, finally,   $L_6$ and $L_8$  are like
$K_2$ and $K_4$ respectively.  Putting together all these facts we obtain \eqref{h3O}.

Similarly to the case of lower derivatives, equation \eqref{fddO} yields
\begin{equation*}\label{h3Ot}
\|\Omega\|_{H^3}\leq P(\|X\|^2_{4}+\|F(X)\|^2_{L^\infty}+\||N|^{-1}\|_{L^\infty})\|\Omega\|_{H^2}.
\end{equation*}
To finish  it remains to show the corresponding inequality for derivatives of fourth order:
\begin{equation}\label{h4O}
\|\Omega\|_{H^4}\leq P(\|X\|^2_{4}+\|F(X)\|^2_{L^\infty}+\||N|^{-1}\|_{L^\infty})\|\Omega\|_{H^3}.
\end{equation}
Identity \eqref{fduO} allows us to point out the most singular terms in $\dau^4\Omega$:
$$
M_1=\frac{A_\mu}{2\pi}PV\int_{\R^2}
\frac{X(\al)-X(\al-\beta)}{|X(\al)-X(\al-\beta)|^3}\wedge\omega(\al-\beta)d\beta \cdot\dau^4 X(\al),
$$
$$
M_2=\frac{A_\mu}{2\pi}PV\int_{\R^2}
\frac{\dau^3 X(\al)-\dau^3 X(\al-\beta)}{|X(\al)-X(\al-\beta)|^3}\wedge\omega(\al-\beta)d\beta \cdot\dau X(\al),
$$
$$
M_3=-\frac{3A_\mu}{4\pi}PV\int_{\R^2}
C(\al,\beta)\frac{X(\al)- X(\al-\beta)}{|X(\al)-X(\al-\beta)|^5}\wedge\omega(\al-\beta)d\beta \cdot\dau X(\al),
$$
with $C(\al,\beta)=(X(\al)- X(\al-\beta))\cdot(\dau^3 X(\al)-\dau^3 X(\al-\beta))$, and
$$
M_4=\frac{A_\mu}{2\pi}PV\int_{\R^2}
\frac{X(\al)-X(\al-\beta)}{|X(\al)-X(\al-\beta)|^3}\wedge\dau^3\omega(\al-\beta)d\beta\cdot\dau X(\al).
$$
Then in order to estimate $M_1$ we start with $\|M_1\|_{L^2}\leq C K \|\dau^4 X\|_{L^2}$
where$$K=\sup_{\al}\Big|PV\int_{\R^2}\frac{X(\al)-X(\al-\beta)}{|X(\al)
-X(\al-\beta)|^3}\wedge\omega(\al-\beta)d\beta \Big|.
$$
Following ref. \cite{DP} we have:
$$
K\leq O_1+O_2+O_3+O_4+O_5
$$
where
$$
O_1=\sup_{\al}\Big|PV\int_{|\beta|>1}\frac{X(\al)-X(\al-\beta)}{|X(\al)
-X(\al-\beta)|^3}\wedge\omega(\al-\beta)d\beta \Big|,
$$$$
O_2=\sup_{\al}\Big|\int_{|\beta|<1}\frac{X(\al)-X(\al-\beta)-\grad X(\al)\cdot\beta}{|X(\al)-X(\al-\beta)|^3}\wedge\omega(\al-\beta)d\beta \Big|,
$$$$
O_3=\sup_{\al}\Big|\int_{|\beta|<1}\grad X(\al)\cdot\beta[\frac{1}{|X(\al)
-X(\al-\beta)|^3}-\frac{1}{|\grad X(\al)\cdot\beta|^3}]\wedge\omega(\al-\beta)d\beta \Big|,
$$$$
O_4=\sup_{\al}\Big| \int_{|\beta|<1}\frac{\grad X(\al)\cdot\beta}{|\grad X(\al)\cdot\beta|^3}\wedge(\omega(\al-\beta)-\omega(\al))d\beta \Big|,
$$$$
O_5=\sup_{\al}\Big|PV\int_{|\beta|<1}\frac{\grad X(\al)\cdot\beta}{|\grad X(\al)\cdot\beta|^3}\wedge \omega(\al)d\beta \Big|.
$$
An  integration by parts  in $O_1$ yields
\begin{align*}
O_1\leq &C\|\grad X \|^2_{L^\infty}\|F(X)\|^3_{L^\infty}\sup_{\al}(\int_{|\beta|>1}\frac{|\Omega(\al-\beta)|}{|\beta|^3}d\beta+\int_{|\beta|=1}|\Omega(\al-\beta)|dl(\beta))\\
\leq & C\|\grad X \|^2_{L^\infty}\|F(X)\|^3_{L^\infty}\|\Omega\|_{L^\infty},
\end{align*}
and  Sobolev's embedding allows us to conclude.

Regarding  $O_2$ we have
$$
O_2\leq \|X-(\al,0)\|_{C^{2,\delta}}\|F(X)\|^3_{L^\infty}\|\omega\|_{L^\infty}\big|\int_{|\beta|<1}|\beta|^{2-\delta}d\beta\big|
$$
 and then the estimate, $\|\omega\|_{C^\delta}\leq C\|\omega\|_{H^2}$ for $0<\delta<1$, gives the desired control. Using \eqref{sioe} and after some straightforward algebraic manipulations we get a similar inequality for $O_3$. Next we have
$$
O_4\leq C \|X-(\al,0)\|^4_{C^1}\||N|^{-1}\|^3_{L^\infty}\|\omega\|_{C^\delta}\big|\int_{|\beta|<1}|\beta|^{2-\delta}d\beta\big|,
$$
giving us also the same  estimate. Furthermore it is easy to prove that $O_5=0$.

 Next we consider the term $M_2$ with the splitting: $M_2=Q_1+Q_2+Q_3$ where
$$
Q_1=\frac{A_\mu}{2\pi}\int_{|\beta|>1}
\frac{\dau^3 X(\al)-\dau^3 X(\al-\beta)}{|X(\al)-X(\al-\beta)|^3}
\wedge\omega(\al-\beta)d\beta \cdot\dau X(\al),
$$
$$
Q_2=\frac{A_\mu}{2\pi}\int_{|\beta|<1}
\frac{\dau^3 X(\al)-\dau^3 X(\al-\beta)}{|X(\al)-X(\al-\beta)|^3}
\wedge(\omega(\al-\beta)-\omega(\al))d\beta \cdot\dau X(\al),
$$
$$
Q_3=\frac{A_\mu}{2\pi}PV\int_{|\beta|<1}
\frac{\dau^3 X(\al)-\dau^3 X(\al-\beta)}{|X(\al)-X(\al-\beta)|^3}
d\beta \wedge\omega(\al)\cdot\dau X(\al).
$$
 The term $Q_1$ can be estimated as before, regarding $Q_2$ we can use the identity $$\dau^3 X(\al)-\dau^3 X(\al-\beta)=\int_0^1\grad \dau^3 X(\al+(s-1)\beta)ds\cdot \beta$$ and the control of  $Q_3$ can be approached as we did with the operator in \eqref{operadorgeneral3}. Similarly with $M_3$, whether  $M_4$ is analogous to $J_6$, and all these observations together allow us to obtain  \eqref{h4O}.


\section{Controlling  the Birkhoff-Rott integral}

Here we consider estimates for the Birkhoff-Rott integral along a non-selfintersecting surface. Let us assume that $\grad(X(\al)-(\al,0))\in H^{k}(\R^2)$ for $k\geq 3$, and that both $F(X)$ and $|N|^{-1}$ are in $L^\infty$ where
$$F(X)(\al,\beta)=|\beta|/|X(\al)-X(\al-\beta)|\qquad \mbox{and}\qquad N(\al)=\dau X(\al)\wedge \dad X(\al).$$
The main purpose of this section is to prove the following estimate:
\begin{eqnarray}\label{estimateBR}
\|BR(X,\omega)\|_{H^{k-1}}\leq  P(\|X\|^2_{k}+\|F(X)\|^2_{L^\infty}+\||N|^{-1}\|_{L^\infty}),
\end{eqnarray}
for $k\geq 4$. Here we shall show it when $k=4$, because the other cases, $k>4$, follow by similar arguments. We rewrite $BR$ in the following manner:
$$
BR(X,\omega)(\al,t)=-\frac{1}{4\pi}PV\int_{\R^2}\frac{X(\al)-X(\beta)}{|X(\al)-X(\beta)|^3}\wedge (\dbd(\Omega\dbu X)-\dbu(\Omega\dbd X))(\beta) d\beta,
$$
 which together with the estimates about $\Omega$ in section 5 and also about the operator $\mathcal{T}_1$ in the appendix, yields
\begin{eqnarray*}
\|BR(X,\omega)\|_{L^2}\leq  P(\|X\|^2_{4}+\|F(X)\|^2_{L^\infty}+\||N|^{-1}\|_{L^\infty}).
\end{eqnarray*}
To estimate  derivatives of order 3 we consider  $\dai^3 (BR(X,\omega))$, and observe that the most dangerous terms are given by
 $$
I_1=-\frac{1}{4\pi}PV\int_{\R^2}\frac{(\dai^3X(\al)-\dai^3X(\al-\beta))\wedge \omega(\al-\beta)}{|X(\al)-X(\al-\beta)|^3}d\beta,
$$
$$
I_2=\frac{3}{4\pi}PV\!\!\int_{\R^2}\!(X(\al)\!-\!X(\al\!-\!\beta))\!\wedge\!\omega(\al\!-\!\beta)\frac{(X(\al)\!-\!X(\al\!-\!\beta))\!\cdot\!(\dai^3X(\al)\!-\!\dai^3X(\al\!-\!\beta))}{|X(\al)\!-\!X(\al\!-\!\beta)|^5}d\beta
$$
$$
I_3=-\frac{1}{4\pi}PV\int_{\R^2}\frac{(X(\al)-X(\al-\beta))\wedge (\dai^3\omega)(\al-\beta)}{|X(\al)-X(\al-\beta)|^3} d\beta.
$$
 In the appendix we find all the ingredients needed to estimate these terms $I_j$ while the remainder in $\dai^3 (BR(X,\omega))$ is easily bounded, namely:  in $I_3$ we can recognize an operator with the form of $\mathcal{T}_1$ in \eqref{operadorgeneral}, so the estimate for $\omega$ in section 5 gives the desired control for $I_3$.
 Regarding  $I_1$ we may use the splitting $I_1=J_1+J_2$ where
$$
J_1=\frac{1}{4\pi}\int_{\R^2}\frac{(\dai^3X(\al)-\dai^3X(\al-\beta))\wedge (\omega(\al)-\omega(\al-\beta))}{|X(\al)-X(\al-\beta)|^3}d\beta,
$$$$
J_2=\frac{\omega(\al)}{4\pi}\wedge PV\int_{\R^2}\frac{(\dai^3X(\al)-\dai^3X(\al-\beta)) }{|X(\al)-X(\al-\beta)|^3}d\beta.
$$
  Then the identity $\dai^3X(\al)-\dai^3X(\al-\beta)=\beta\cdot \int_0^1\grad\dai^3 X(\al+(s-1)\beta)ds$ allows us to find in $J_1$ a kernel of degree $-1$ which we know how to handle (see appendix). One  use the estimate for $\mathcal{T}_3$ \eqref{operadorgeneral3} to deal with $J_2$ and we proceed similarly to control $I_2$.

\section{In search of the Rayleigh-Taylor condition}

As it was pointed out in section 4 (outline of the proof) our approach is based on energy estimates and a crucial step is to characterize those terms involving higher derivatives which are controlled because they have the appropriated sign. In our terminology they constitute the Rayleigh-Taylor condition, which is supposed to holds at time $T=0$, being an important part of the proof to show that it prevails under the evolution.

 Let us introduce the notation
$$
|||X|||_k^2=\|X\|^2_k+\|F(X)\|^2_{L^\infty}+\||N|^{-1}\|_{L^\infty}
$$
where
\begin{equation}\label{nl3l2s}
\|X\|_k=\|X_1-\al_1\|_{L^3}+\|X_2-\al_2\|_{L^3}+\|X_3\|_{L^2}+\|\grad(X-(\al,0))\|^2_{H^{k-1}},
\end{equation}
and
$$
\|\grad(X-(\al,0))\|^2_{H^{k-1}}=\|\grad(X-(\al,0))\|^2_{L^2}+\|\dau^k(X-(\al,0))\|^2_{L^2}+\|\dad^k(X-(\al,0))\|^2_{L^2}.
$$

In order to justify the formula
\begin{align*}
\frac{d}{dt}\|X\|^2_{k}(t)&\leq
\D-\sum_{i=1,2}\frac{2^{3/2}}{(\mu_1\!+\!\mu_2)}\,\int_{\R^2} \frac{\sigma(\al,t)}{|\grad X(\al,t)|^{3}}
 \dai^{k} X(\al,t)\cdot \la(\dai^{k} X)(\al,t) d\al\\
&\quad + P(|||X|||_{k}(t)),
\end{align*}
(here $k\geq 4$, although for the sake of simplicity we will present the explicit computations when $k=4$, leaving the other cases as an exercise for the interested reader), it will be convenient to make use of the following tools, giving us different kind of cancelations, and which constitute our particular  bestiary of formulas for this paper:

  From the definition of the isothermal parameterization we have the identities:
\begin{equation}\label{b1}
|\dau X|^2=|\dad X|^2,
\end{equation}
\begin{equation}\label{b2}
\dau X\cdot\dad X=0,
\end{equation}
which yield
\begin{equation}\label{b3}
\frac12\Delta(|\dau X|^2)=|\dau\dad X|^2-\dau^2X\cdot\dad^2X,
\end{equation}
\begin{align}
\begin{split}\label{b4a}
\dau^4X\cdot\dau X&=-3\dau^3X\cdot\dau^2X+(\dau^2\Delta^{-1}\dau)(|\dau\dad X|^2-\dau^2X\cdot\dad^2X),
\end{split}
\end{align}
\begin{align}
\begin{split}\label{b4b}
\dad^4X\cdot\dad X&=-3\dad^3X\cdot\dad^2X+(\dad^2\Delta^{-1}\dad)(|\dau\dad X|^2-\dau^2X\cdot\dad^2X).
\end{split}
\end{align}
Using \eqref{b2} and \eqref{b3} we obtain:
\begin{align}
\begin{split}\label{b4c}
\dau^4X\cdot\dad X&=-2\dau^3X\cdot\dau\dad X-\dau^2\dad X\cdot\dau^2 X\\
&\quad -(\dau\dad\Delta^{-1}\dau)(|\dau\dad X|^2-\dau^2X\cdot\dad^2X),
\end{split}
\end{align}
\begin{align}
\begin{split}\label{b4d}
\dad^4X\cdot\dau X&=-2\dad^3X\cdot\dau\dad X-\dad^2\dau X\cdot\dad^2 X\\
&\quad -(\dau\dad\Delta^{-1}\dad)(|\dau\dad X|^2-\dau^2X\cdot\dad^2X).
\end{split}
\end{align}
 And Sobolev inequalities imply  that if $\grad(X-(\al,0))\in H^3$ then $\dai^4X\cdot\daj X\in H^3$ for $i,j=1,2$.

With the help of the estimates above we may now afford the task of determining $\sigma$. There is a part that  may be  considered as a mere ``algebraic'' manipulation to detect the relevant characters and, in so doing, we disregard many terms because they are of lower order in the sense of Sobolev spaces. At the end, we shall present how to deal with those lower order terms, if not for the whole collection of them, at least for the ones that we may consider to be the most ``dangerous'' characters. Here it is convenient to recommend the reader  our previous works \cite{DP,ADP} where similar estimates were carried out.

\subsection{Low order norms}

  Since $X_i(\al)\to \al_i$ for $i=1,2$ at infinity, let us consider the evolution of the $L^3$ norm. That is
\begin{align*}
\frac13\frac{d}{dt}\|X_1-\al_1\|_{L^3}^3(t)&=\int_{\R^2}|X_1-\al_1|(X_1-\al_1)X_{1t}d\al=I_1+I_2+I_3,
\end{align*}
where
$$
I_1=\int_{\R^2}|X_1-\al_1|(X_1-\al_1)BR_1d\al,
$$
$$
I_2=\int_{\R^2}|X_1-\al_1|(X_1-\al_1)C_1\dau X_1d\al,\quad I_3=\int_{\R^2}|X_1-\al_1|(X_1-\al_1)C_2\dad X_1d\al.
$$
Then we have
$$
I_1\leq \|X_1-\al_1\|^2_{L^3}\|BR\|_{L^3}\leq C(\|X_1-\al_1\|^3_{L^3}+\|BR\|_{L^\infty}\|BR\|^2_{L^2}),
$$
and Sobolev estimates, together with \eqref{estimateBR}, yield the appropriate control in terms of $P(|||X|||_k)$.

Next since $\dau X_1\to 1$ as $\al\to\infty$,  we have
$$
I_2\leq \|\dau X_1\|_{L^\infty}\|X_1-\al_1\|^2_{L^3}\|C_1\|_{L^3},
$$
and it remains to get control of $C_1$. Using \eqref{C1} we introduce the splitting $C_1=\sum_{j=1}^4C_1^j$, where

$$
C_1^1(\al)=\frac{1}{2\pi}\int_{\R^2}\frac{\al_1-\beta_1}{|\al-\beta|^2}BR_{\beta_2}\cdot\frac{X_{\beta_2}}{|X_{\beta_2}|^2}d\beta,
\quad
C_1^2(\al)=-\frac{1}{2\pi}\int_{\R^2}\frac{\al_1-\beta_1}{|\al-\beta|^2}BR_{\beta_1}\cdot \frac{X_{\beta_1}}{|X_{\beta_2}|^2}d\beta,
$$
$$
C_1^3(\al)=-\frac{1}{2\pi}\int_{\R^2}\frac{\al_2-\beta_2}{|\al-\beta|^2}BR_{\beta_1}\cdot\frac{X_{\beta_2}}{|X_{\beta_1}|^2}d\beta,
\quad
C_1^4(\al)=-\frac{1}{2\pi}\int_{\R^2}\frac{\al_1-\beta_1}{|\al-\beta|^2}BR_{\beta_2}\cdot \frac{X_{\beta_1}}{|X_{\beta_1}|^2}d\beta.
$$
We shall show how control $C^1_1$, because the estimates for the other terms follow by similar arguments. Integrating by parts one obtain
$C_1^1=D_1+D_2$ where
$$
D_1\!=\!\frac{-1}{2\pi}\int_{\R^2}\frac{\al_1\!-\!\beta_1}{|\al\!-\!\beta|^2}BR\!\cdot\!\dbd\Big(\frac{X_{\beta_2}}{|X_{\beta_2}|^2}\Big)d\beta,
\,\,
D_2\!=\!-\frac{1}{\pi}PV\!\!\int_{\R^2}\!\!\frac{(\al_1\!-\!\beta_1)(\al_2\!-\!\beta_2)}{|\al\!-\!\beta|^4}BR\!\cdot\!\frac{X_{\beta_2}}{|X_{\beta_2}|^2}d\beta.
$$
Regarding $D_1$ we write $D_1=E_1+E_2$ where
$$E_1\!=\!\frac{-1}{2\pi}\int_{|\beta|<1}\frac{\beta_1}{|\beta|^2}BR(\al\!-\!\beta)\!\cdot\!\dbd\Big(\frac{X_{\beta_2}}{|X_{\beta_2}|^2}\Big)(\al\!-\!\beta)d\beta,$$
$$E_2\!=\!\frac{-1}{2\pi}\int_{|\beta|>1}\frac{\beta_1}{|\beta|^2}BR(\al\!-\!\beta)\!\cdot\!\dbd\Big(\frac{X_{\beta_2}}{|X_{\beta_2}|^2}\Big)(\al\!-\!\beta)d\beta.
$$
Then Minkowski and Young inequalities yield respectively
$$\|E_1\|_{L^3}\leq C\|BR\!\cdot\!\dbd\Big(\frac{X_{\beta_2}}{|X_{\beta_2}|^2}\Big)\|_{L^3}\leq P(|||X|||_4), $$
$$\|E_2\|_{L^3}\leq C\|BR\!\cdot\!\dbd\Big(\frac{X_{\beta_2}}{|X_{\beta_2}|^2}\Big)\|_{L^1}\leq C\|BR\|_{L^2}\|\dbd \Big(\frac{X_{\beta_2}}{|X_{\beta_2}|^2}\Big)\|_{L^2}\leq P(|||X|||_4),$$
and the desired control is achieved. In the term $D_2$ we have a double Riesz transform and the standard Calderon-Zygmund theory yields
$$
\|D_2\|_{L^3}\leq C\|BR\cdot\frac{X_{\beta_2}}{|X_{\beta_2}|^2}\|_{L^3}\leq C\||X_{\beta_2}|^{-1}\|_{L^\infty}\|BR\|_{L^3}\leq P(|||X|||_4).
$$
 The estimate for $I_3$ follows on similar path, and the case of the second coordinate is also identical:
\begin{align*}
\frac13\frac{d}{dt}\|X_2-\al_2\|_{L^3}^3(t)\leq P(|||X|||_4).
\end{align*}
Regarding  the third coordinate we have a stronger decay because of the asymptotic flatness hypothesis:
\begin{align*}
\frac12\frac{d}{dt}\|X_3\|_{L^2}^2(t)&=\int_{\R^2}\!\!X_3 BR_3d\al+\int_{\R^2}\!\!X_3C_1\dau X_3 d\al+\int_{\R^2}\!\!X_3C_2\dad X_3d\al\\
&=\int_{\R^2}\!\!X_3 BR_3d\al-\frac12\int_{\R^2}(\dau C_1+\dad C_2)|X_3|^2 d\al,
\end{align*}
 therefore the use of Sobolev's embedding in the formulas for $C_1$ \eqref{C1} and $C_2$ \eqref{C2}, together with the estimates for $BR$ \eqref{estimateBR}, allows us to obtain:
\begin{align*}
\frac12\frac{d}{dt}\|X_3\|_{L^2}^2(t)\leq P(|||X|||_4).
\end{align*}

Once we have control of higher order derivatives, we can use the estimates of the appendix to get
\begin{equation*}
\frac12\frac{d}{dt}\|\grad(X-(\al,0))\|_{L^2}^2(t)\leq P(|||X|||_4).
\end{equation*}

\subsection{Higher order norms}
Let us  consider now
\begin{align}
\begin{split}\label{dauX}
\frac12\frac{d}{dt}\|\dau^4X\|_{L^2}^2(t)&=\int_{\R^2}\!\!\dau^4X\cdot\dau^4BR(X,\omega)d\al\\
&\quad+\int_{\R^2}\!\!\dau^4X\cdot\dau^4(C_1\dau X)d\al+\int_{\R^2}\!\!\dau^4X\cdot\dau^4(C_2\dad X)d\al\\
&=I_1+I_2+I_3,
\end{split}
\end{align}
The higher order terms in $I_2$ and $I_3$ are given by
$$
J_1=\int_{\R^2}\!\!C_1\dau^4X\cdot\dau^5 Xd\al,\qquad J_2=\int_{\R^2}\!\!\dau^4X\cdot\dau X  \dau^4C_1d\al
$$
$$
J_3=\int_{\R^2}\!\!C_2\dau^4X\cdot\dau^4\dad Xd\al,\qquad J_4=\int_{\R^2}\!\!\dau^4X\cdot\dad X  \dau^4C_2d\al
$$
Integration by parts yields
$$
J_1+J_3=-\frac{1}{2}\int_{\R^2}\!\!(\dau C_1+\dad C_2)|\dau^4X|^2d\al
$$
and therefore
$$
J_1+J_3\leq \frac12(\|\dau C_1\|_{L^\infty}+\|\dad C_2\|_{L^\infty})\|\dau^4X\|_{L^2}^2\leq  P(|||X|||_{4}).
$$
Then in $J_2$ we use \eqref{b4a} to get
$$
J_2=-\int_{\R^2}\!\!\dau(\dau^4X\cdot\dau X)  \dau^3C_1d\al\leq \|\dau(\dau^4X\cdot\dau X)\|_{L^2}\|\dau^3C_1\|_{L^2}.
$$
Whether in $J_4$ we use \eqref{b4c} to obtain
$$
J_4=-\int_{\R^2}\!\!\dau(\dau^4X\cdot\dad X)  \dau^3C_2d\al\leq \|\dau(\dau^4X\cdot\dad X)\|_{L^2}\|\dau^3C_2\|_{L^2}.
$$
 From formulas \eqref{C1},\eqref{C2} one  realizes that $C_1$ and $C_2$ are at the same level than Birkhoff-Rott \eqref{BR}, and, therefore,  we can use the estimates for $BR$ \eqref{estimateBR} to control $\|\dau^3C_i\|_{L^2}$, $i=1,2$. Then formulas \eqref{b4a} and \eqref{b4c} indicate how to estimate $\|\dau(\dau^4X\cdot\dai X)\|_{L^2}$, $i=1,2$. That is we have:
$$
J_2+J_4\leq P(|||X|||_{4}).
$$
In $I_1$ the most singular terms are given by
$$
J_5=-\frac{1}{4\pi}PV\int_{\R^2}\int_{\R^2}\dau^4X(\al)\cdot\frac{(\dau^4X(\al)-\dau^4X(\beta))\wedge \omega(\beta)}{|X(\al)-X(\beta)|^3} d\al d\beta,
$$
\begin{equation}\label{J6}
J_6=\frac{3}{4\pi}PV\!\!\int_{\R^2}\!\!\int_{\R^2}\!\!\dau^4X(\al)\!\cdot\!(X(\al)\!-\!X(\beta))\!\wedge\!\omega(\beta)\frac{(X(\al)\!-\!X(\beta))\!\cdot\!(\dau^4X(\al)\!-\!\dau^4X(\beta))}{|X(\al)-X(\beta)|^5} d\al d\beta
\end{equation}
$$
J_7=-\frac{1}{4\pi}PV\int_{\R^2}\int_{\R^2}\dau^4X(\al)\cdot
\frac{(X(\al)-X(\beta))\wedge (\dau^4\omega)(\beta)}{|X(\al)-X(\beta)|^3} d\al d\beta.
$$
Let us consider now the splitting $J_5=K_1+K_2$
$$
K_1=-\frac{1}{8\pi}PV\int_{\R^2}\int_{\R^2}\dau^4X(\al)\wedge (\dau^4X(\al)-\dau^4X(\beta))\cdot\frac{ \omega(\beta)+\omega(\al)}{|X(\al)-X(\beta)|^3} d\al d\beta,
$$
$$
K_2=\frac{1}{8\pi}PV\int_{\R^2}\int_{\R^2}\dau^4X(\al)\wedge (\dau^4X(\al)-\dau^4X(\beta))\cdot\frac{\omega(\al)-\omega(\beta)}{|X(\al)-X(\beta)|^3} d\al d\beta,
$$
Next we exchange $\al$ and $\beta$ in $K_1$ to get
\begin{align*}
K_1&=\frac{1}{8\pi}PV\int_{\R^2}\int_{\R^2}\dau^4X(\beta)\wedge (\dau^4X(\al)-\dau^4X(\beta))\cdot\frac{ \omega(\beta)+\omega(\al)}{|X(\al)-X(\beta)|^3} d\al d\beta\\
&=\frac{-1}{16\pi}PV\int_{\R^2}\int_{\R^2}(\dau^4X(\al)-\dau^4X(\beta))\wedge (\dau^4X(\al)-\dau^4X(\beta))\cdot\frac{ \omega(\beta)+\omega(\al)}{|X(\al)-X(\beta)|^3} d\al d\beta
\end{align*}
and therefore we can conclude that $K_1=0$. In $K_2$ we find a singular integral with a kernel of degree $-2$
$$
K_2=-\frac{1}{8\pi}PV\int_{\R^2}\dau^4X(\al)\cdot \int_{\R^2}\dau^4X(\beta)\wedge\frac{\omega(\al)-\omega(\beta)}{|X(\al)-X(\beta)|^3} d\beta d\al,
$$
and as it is proved in the appendix we have
$$
K_2\leq P(|||X|||_{4}).
$$
Let us now  decompose $J_6=K_3+K^1_4+K_4^2+K_5^1+K_5^2$ where
$$
K_3=\frac{3}{4\pi}PV\int_{\R^2}\int_{\R^2}\!\!\dau^4X(\al)\cdot(X(\al)-X(\beta))\wedge\omega(\beta)\frac{A(\al,\beta)\cdot(\dau^4X(\al)\!-\!\dau^4X(\beta))}{|X(\al)-X(\beta)|^5} d\al d\beta
$$
with $A(\al,\beta)=X(\al)-X(\beta)-\grad X(\al)(\al-\beta)$,
$$
K^i_4=\frac{-3}{4\pi}PV\int_{\R^2}\int_{\R^2}\!\!\dau^4X(\al)\cdot(X(\al)-X(\beta))\wedge\omega(\beta)\frac{(\al_i-\beta_i)(\dai X(\al)\!-\!\dai X(\beta))\cdot \dau^4X(\beta)}{|X(\al)-X(\beta)|^5} d\al d\beta
$$
$$
K^i_5=\frac{3}{4\pi}PV\int_{\R^2}\int_{\R^2}\!\!\dau^4X(\al)\cdot(X(\al)-X(\beta))\wedge\omega(\beta)\frac{(\al_i\!-\!\beta_i)(\dai X(\al)\cdot\dau^4X(\al)\!-\!\dai X(\beta)\cdot\dau^4X(\beta))}{|X(\al)-X(\beta)|^5} d\al d\beta
$$
In $K_3$ and $K_4^i$ we find kernels of degree $-2$ and, as it is shown in the appendix, they behave as a Riesz transform acting on $\dau^4X$. In $K_5^i$ the kernels have degree $-3$ and act as a $\Lambda$ operator  on $\dai X\cdot\dau^4X$. Then using the formulas \eqref{b4a} and \eqref{b4c} we get finally the desired estimate.

We will find the R-T condition in $J_7$. Let us take $J_7=K_6+K_7$ where
$$
K_6=-\frac{1}{4\pi}PV\int_{\R^2}\dau^4X(\al)\cdot\int_{\R^2}
\big(\frac{(X(\al)-X(\beta))}{|X(\al)-X(\beta)|^3}-\frac{\grad X(\al)(\al-\beta)}{|\grad X(\al)(\al-\beta)|^3}\big) \wedge (\dau^4\omega)(\beta)d\beta d\al,
$$
$$
K_7=-\frac{1}{4\pi}PV\int_{\R^2}\dau^4X(\al)\cdot\int_{\R^2}
\frac{\grad X(\al)(\al-\beta)}{|\grad X(\al)(\al-\beta)|^3}\wedge (\dau^4\omega)(\beta)d\beta d\al.
$$
The term $K_6$ is controlled by \eqref{operadorregular} in the appendix. Using \eqref{b1} and \eqref{b2} we get
$$
K_7=-\frac12PV\int_{\R^2}\frac{\dau^4X(\al)}{|\dau X(\al)|^3}\cdot (\dau X(\al)\wedge R_1(\dau^4\omega)(\al)+\dad X(\al)\wedge R_2(\dau^4\omega)(\al)) d\al.
$$
Formula \eqref{oOf} help us to detect the most singular terms inside  $K_7$, which will be denoted by $L_i$, $i=1,...,8$ and are the following:
$$
L_1=-\frac12PV\int_{\R^2}\dau^4X(\al)\cdot \frac{\dau X(\al)}{|\dau X(\al)|^3}\wedge R_1(\dau^4\dad\Omega \dau X)(\al)d\al,
$$
$$
L_2=-\frac12PV\int_{\R^2}\dau^4X(\al)\cdot \frac{\dau X(\al)}{|\dau X(\al)|^3}\wedge R_1(\dad\Omega \dau^5 X)(\al)d\al,
$$
$$
L_3=\frac12PV\int_{\R^2}\dau^4X(\al)\cdot \frac{\dau X(\al)}{|\dau X(\al)|^3}\wedge R_1(\dau^5\Omega \dad X)(\al)d\al,
$$$$
L_4=\frac12PV\int_{\R^2}\dau^4X(\al)\cdot \frac{\dau X(\al)}{|\dau X(\al)|^3}\wedge R_1(\dau\Omega \dau^4\dad X)(\al)d\al,
$$$$
L_5=-\frac12PV\int_{\R^2}\dau^4X(\al)\cdot \frac{\dad X(\al)}{|\dad X(\al)|^3}\wedge R_2(\dau^4\dad\Omega \dau X)(\al)d\al,
$$$$
L_6=-\frac12PV\int_{\R^2}\dau^4X(\al)\cdot \frac{\dad X(\al)}{|\dad X(\al)|^3}\wedge R_2(\dad\Omega \dau^5 X)(\al)d\al,
$$$$
L_7=\frac12PV\int_{\R^2}\dau^4X(\al)\cdot \frac{\dad X(\al)}{|\dad X(\al)|^3}\wedge R_2(\dau^5\Omega \dad X)(\al)d\al,
$$$$
L_8=\frac12PV\int_{\R^2}\dau^4X(\al)\cdot \frac{\dad X(\al)}{|\dad X(\al)|^3}\wedge R_2(\dau\Omega \dau^4\dad X)(\al)d\al.
$$
In $L_1$ we get a kernel of degree $-1$ of the form
$$
L_1=\frac12PV\int_{\R^2}\dau^4X(\al)\cdot \int_{\R^2}\frac{\al_1-\beta_1}{|\al-\beta|^3} \frac{\dau X(\al)}{|\dau X(\al)|^3}\wedge (\dau X(\al)-\dau X(\beta))\dau^4\dad\Omega(\beta) d\beta d\al,
$$
 which can be estimated integrating by parts throughout  $\dau^4\dad\Omega$; also the term  $L_7$ follows in a similar manner. In order to estimate $L_2,$ $L_4,$ $L_6$ and $L_8$ we realize that they can be written like \eqref{commomega} in the appendix plus commutators of the form \eqref{commderR}. Next we have to deal  with $L_3$ and $L_5$:
 With $L_3$ we proceed as follows
$$L_3\leq \widetilde{L}_3+\||\dau X|^{-2}\|_{L^\infty}\|\dau^4X\|_{L^2}\|R_1(\dau^5\Omega \dad X)-R_1(\dau^5\Omega)\dad X\|_{L^2}
$$
where $\widetilde{L}_3$ is given by
\begin{equation}\label{Lt3}
\widetilde{L}_3=\frac12PV\int_{\R^2}\dau^4X(\al)\cdot \frac{N(\al)}{|\dau X(\al)|^3}(R_1\dau)(\dau^4\Omega)(\al)d\al,
\end{equation}
and the commutator estimates \eqref{commderR} show that it only remains to control $\widetilde{L}_3$. We use now formula \eqref{fduO} to get $\widetilde{L}_3=M_1+M_2$ where
$$
M_1=-A_\rho PV\int_{\R^2}\dau^4X(\al)\cdot \frac{N(\al)}{|\dau X(\al)|^3}(R_1\dau)(\dau^4 X_3)(\al)d\al,
$$
and
$$
M_2=-A_\mu PV\int_{\R^2}\dau^4X(\al)\cdot \frac{N(\al)}{|\dau X(\al)|^3}(R_1\dau)(\dau^3 (BR(X,\omega)\cdot\dau X))(\al)d\al.
$$
Then we write $M_1=O_1+O_2+O_3$ where
$$
O_1=-A_\rho PV\int_{\R^2}\frac{\dau^4X_1}{|\dau X|^3}(\dau X_2\dad X_3-\dau X_3\dad X_2)(R_1\dau)(\dau^4 X_3)d\al,
$$
$$
O_2=-A_\rho PV\int_{\R^2}\frac{\dau^4X_2}{|\dau X|^3}(\dau X_3\dad X_1-\dau X_1\dad X_3)(R_1\dau)(\dau^4 X_3)d\al,
$$
\begin{equation}\label{O3}
O_3=-A_\rho PV\int_{\R^2}\frac{N_3}{|\dau X|^3}\dau^4X_3(R_1\dau)(\dau^4 X_3)d\al.
\end{equation}
Next we consider $O_1=P_1+P_2+P_3$ with
$$
P_1=-A_\rho PV\int_{\R^2}\frac{\dau^4X_1}{|\dau X|^3}\dau X_2(R_1\dau)(\dad X_3\dau^4 X_3)d\al,
$$
$$
P_2=A_\rho PV\int_{\R^2}\frac{\dau^4X_1}{|\dau X|^3}\dad X_2(R_1\dau)(\dau X_3\dau^4 X_3)d\al,
$$
\begin{align*}
P_3=&A_\rho PV\int_{\R^2}\frac{\dau^4X_1}{|\dau X|^3}\dau X_2[(R_1\dau)(\dad X_3\dau^4 X_3)-\dad X_3(R_1\dau)(\dau^4 X_3)]d\al\\
&+A_\rho PV\int_{\R^2}\frac{\dau^4X_1}{|\dau X|^3}\dad X_2
[\dau X_3(R_1\dau)(\dau^4 X_3)-(R_1\dau)(\dau X_3\dau^4 X_3)]d\al
\end{align*}
and the commutator estimate allows us to control the term $P_3$.

Now we use \eqref{b4c} to write $P_1=Q_1+Q_2+Q_3$
$$
Q_1=A_\rho PV\int_{\R^2}\frac{\dau^4X_1}{|\dau X|^3}\dau X_2(R_1\dau)(\dad X_1\dau^4 X_1)d\al,
$$
$$
Q_2=A_\rho PV\int_{\R^2}\frac{\dau^4X_1}{|\dau X|^3}\dau X_2(R_1\dau)(\dad X_2\dau^4 X_2)d\al,
$$
$$
Q_3=A_\rho PV\int_{\R^2}\frac{\dau^4X_1}{|\dau X|^3}\dau X_2(R_1\dau)(\mbox{l.o.t.})d\al.
$$
The term $Q_3$ is easily estimated. Regarding $P_2$ equality  \eqref{b4a} allows us to write $P_2=Q_4+Q_5+Q_6$ where
$$
Q_4=-A_\rho PV\int_{\R^2}\frac{\dau^4X_1}{|\dau X|^3}\dad X_2(R_1\dau)(\dau X_1\dau^4 X_1)d\al,
$$
$$
Q_5=-A_\rho PV\int_{\R^2}\frac{\dau^4X_1}{|\dau X|^3}\dad X_2(R_1\dau)(\dau X_2\dau^4 X_2)d\al,
$$
$$
Q_6=-A_\rho PV\int_{\R^2}\frac{\dau^4X_1}{|\dau X|^3}\dad X_2(R_1\dau)(\mbox{l.o.t.})d\al.
$$
Let us  recall the identity $P_1+P_2=(Q_4+Q_1)+(Q_2+Q_5)+(Q_3+Q_6)$ where $Q_3$ and $
Q_6$ are easily  estimated. With respect to $Q_2+Q_5$ we have
\begin{align*}
Q_2+Q_5&=A_\rho PV\int_{\R^2}\frac{\dau^4X_1}{|\dau X|^3}\dau X_2[(R_1\dau)(\dad X_2\dau^4 X_2)-\dad X_2(R_1\dau)(\dau^4 X_2)]d\al\\
&\quad +A_\rho PV\int_{\R^2}\frac{\dau^4X_1}{|\dau X|^3}\dad X_2[\dau X_2(R_1\dau)(\dau^4 X_2)-(R_1\dau)(\dau X_2\dau^4 X_2)]d\al
\end{align*}
and again the commutator estimates yields the desired control.

Next we have
\begin{align*}
Q_4+Q_1&=A_\rho PV\int_{\R^2}\frac{\dau^4X_1}{|\dau X|^3}\dad X_2[\dau X_1(R_1\dau)(\dau^4 X_1)-(R_1\dau)(\dau X_1\dau^4 X_1)]d\al\\
&\quad + A_\rho PV\int_{\R^2}\frac{\dau^4X_1}{|\dau X|^3}\dau X_2[(R_1\dau)(\dad X_1\dau^4 X_1)-\dad X_1(R_1\dau)(\dau^4 X_1)]d\al\\
&\quad - A_\rho PV\int_{\R^2}\frac{N_3}{|\dau X|^3}\dau^4X_1(R_1\dau)(\dau^4 X_1)d\al.
\end{align*}
The first two integrals above are easily handled allowing us to get
\begin{align}\label{O1}
\begin{split}
O_1=P_1+P_2+P_3&\leq P(|||X|||_{4})-A_\rho PV\int_{\R^2}\frac{N_3}{|\dau X|^3}\dau^4X_1(R_1\dau)(\dau^4 X_1)d\al.
\end{split}
\end{align}
For the term $O_2$ we proceed in a similar manner, first we check that $O_2=P_4+P_5+P_6$
$$
P_4=A_\rho PV\int_{\R^2}\frac{\dau^4X_2}{|\dau X|^3}\dau X_1(R_1\dau)(\dad X_3\dau^4 X_3)d\al,
$$
$$
P_5=-A_\rho PV\int_{\R^2}\frac{\dau^4X_2}{|\dau X|^3}\dad X_1(R_1\dau)(\dau X_3\dau^4 X_3)d\al,
$$
\begin{align*}
P_6=&A_\rho PV\int_{\R^2}\frac{\dau^4X_2}{|\dau X|^3}\dau X_1[\dad X_3(R_1\dau)(\dau^4 X_3)-(R_1\dau)(\dad X_3\dau^4 X_3)]d\al\\
&+A_\rho PV\int_{\R^2}\frac{\dau^4X_2}{|\dau X|^3}\dad X_1[(R_1\dau)(\dau X_3\dau^4 X_3)-\dau X_3(R_1\dau)(\dau^4 X_3)]d\al.
\end{align*}
We control $P_6$ as before. Regarding $P_4$ we use \eqref{b4c} to write it in the form $P_4=S_1+S_2+S_3$ where:
$$
S_1=-A_\rho PV\int_{\R^2}\frac{\dau^4X_2}{|\dau X|^3}\dau X_1(R_1\dau)(\dad X_1\dau^4 X_1)d\al,
$$
$$
S_2=-A_\rho PV\int_{\R^2}\frac{\dau^4X_2}{|\dau X|^3}\dau X_1(R_1\dau)(\dad X_2\dau^4 X_2)d\al,
$$
$$
S_3=-A_\rho PV\int_{\R^2}\frac{\dau^4X_2}{|\dau X|^3}\dau X_1(R_1\dau)(\mbox{l.o.t.})d\al.
$$
The identity \eqref{b4a} allows us to write $P_5=S_4+S_5+S_6$ where:
$$
S_4=A_\rho PV\int_{\R^2}\frac{\dau^4X_2}{|\dau X|^3}\dad X_1(R_1\dau)(\dau X_1\dau^4 X_1)d\al,
$$
$$
S_5=A_\rho PV\int_{\R^2}\frac{\dau^4X_2}{|\dau X|^3}\dad X_1(R_1\dau)(\dau X_2\dau^4 X_2)d\al,
$$$$
S_6=A_\rho PV\int_{\R^2}\frac{\dau^4X_2}{|\dau X|^3}\dad X_1(R_1\dau)(\mbox{l.o.t.})d\al.
$$
Next, we reorganize the sum in the form $P_4+P_6=(S_1+S_4)+(S_2+S_5)+(S_3+S_6)$ where the term $S_3+S_6$ can be easily estimated. Regarding $S_1+S_4$ we have
\begin{align*}
S_1+S_4&=A_\rho PV\int_{\R^2}\frac{\dau^4X_2}{|\dau X|^3}\dau X_1[\dad X_1(R_1\dau)(\dau^4 X_1)-(R_1\dau)(\dad X_1\dau^4 X_1)]d\al\\
&\quad+A_\rho PV\int_{\R^2}\frac{\dau^4X_2}{|\dau X|^3}\dad X_1[(R_1\dau)(\dau X_1\dau^4 X_1)-\dau X_1(R_1\dau)(\dau^4 X_1)]d\al
\end{align*}
and the commutator estimates gives us precise control.

Let us consider now
\begin{align*}
S_2+S_5&=A_\rho PV\int_{\R^2}\frac{\dau^4X_2}{|\dau X|^3}\dau X_1
[\dad X_2(R_1\dau)(\dau^4 X_2)-(R_1\dau)(\dad X_2\dau^4 X_2)]d\al\\
&\quad+ A_\rho PV\int_{\R^2}\frac{\dau^4X_2}{|\dau X|^3}\dad X_1[(R_1\dau)(\dau X_2\dau^4 X_2)-\dau X_2(R_1\dau)(\dau^4 X_2)]d\al\\
&\quad-A_\rho PV\int_{\R^2}\frac{N_3}{|\dau X|^3}\dau^4X_2(R_1\dau)(\dau^4 X_2)d\al.
\end{align*}
Here again  the commutator estimate  control the first two integrals above, allowing us to conclude that
\begin{align}\label{O2}
\begin{split}
O_2=P_4+P_5+P_6&\leq P(|||X|||_{4})-A_\rho PV\int_{\R^2}\frac{N_3}{|\dau X|^3}\dau^4X_2(R_1\dau)(\dau^4 X_2)d\al.
\end{split}
\end{align}
Furthermore, inequalities \eqref{O1}, \eqref{O2} and \eqref{O3} yield
\begin{align}\label{densidadRdau}
\begin{split}
M_1=O_1+O_2+O_3&\leq P(|||X|||_{4})-A_\rho PV\int_{\R^2}\frac{N_3}{|\dau X|^3}\dau^4X\cdot(R_1\dau)(\dau^4 X)d\al,
\end{split}
\end{align}
and at this point we begin to recognize the Rayleigh-Taylor condition in the non-integrable terms. Let us  return now to the term $M_2$ which can be written in the form
\begin{equation}\label{M22}
M_2= A_\mu PV\int_{\R^2}R_1\big( \frac{\dau^4X\cdot N}{|\dau X|^3}\big)\dau^4 (BR(X,\omega)\cdot\dau X))d\al,
\end{equation}
and whose most dangerous components are given by
$$
O_4=-\frac{A_\mu}{4\pi}PV\int_{\R^2}R_1\big( \frac{\dau^4X\cdot N}{|\dau X|^3}\big)(\al)\int_{\R^2}\frac{\dau^4X(\al)-\dau^4X(\beta)}{|X(\al)-X(\beta)|^3}\wedge \omega(\beta)\cdot\dau X(\al)d\al,
$$
$$
O_5=\frac{3A_\mu}{4\pi}PV\int_{\R^2}R_1\big( \frac{\dau^4X\cdot N}{|\dau X|^3}\big)(\al)\int_{\R^2}B(\al,\beta)(X(\al)-X(\beta))\wedge \omega(\beta)\cdot\dau X(\al)d\al,
$$
with
$$
B(\al,\beta)=\frac{(X(\al)-X(\beta))\cdot(\dau^4X(\al)-\dad^4X(\beta))}{|X(\al)-X(\beta)|^5},
$$
$$
O_6=-\frac{A_\mu}{4\pi}PV\int_{\R^2}R_1\big( \frac{\dau^4X\cdot N}{|\dau X|^3}\big)(\al)\int_{\R^2}\frac{X(\al)-X(\beta)}{|X(\al)-X(\beta)|^3}\wedge \dau^4\omega(\beta)\cdot\dau X(\al)d\al,
$$
and

$$
O_7= A_\mu PV\int_{\R^2}R_1\big( \frac{\dau^4X\cdot N}{|\dau X|^3}\big)(\al)\dau (BR(X,\omega)\cdot \dau^4 X)(\al)d\al.
$$

The remainder terms are less singular and can be estimated with the same methods used before. To deal with $O_4$ we decompose it further $O_4=P_7+P_8$:
$$
P_7=\frac{A_\mu}{4\pi}PV\!\!\int_{\R^2}R_1\big( \frac{\dau^4X\cdot N}{|\dau X|^3}\big)(\al)\int_{\R^2}\frac{\dau^4X(\al)-\dau^4X(\beta)}{|X(\al)-X(\beta)|^3}\cdot \omega(\beta)\wedge(\dau X(\beta)-\dau X(\al))d\beta d\al,
$$
$$
P_8=\frac{A_\mu}{4\pi}PV\!\!\int_{\R^2}R_1\big( \frac{\dau^4X\cdot N}{|\dau X|^3}\big)(\al)\int_{\R^2}\frac{\dau^4X(\al)-\dau^4X(\beta)}{|X(\al)-X(\beta)|^3}\cdot  N(\beta)\dau\Omega(\beta)d\beta d\al,
$$
where in $P_8$ we have used  formula \eqref{oOf} to get $\omega\wedge \dau X=N\dau\Omega$. In the integral (with respect to $\beta$) of $P_7$ we have a kernel of degree $-2$ applied to $4$ derivatives, which   can be estimated easily. Next let us consider $P_8=Q_7+Q_8+Q_9$ where
$$
Q_7=-\frac{A_\mu}{4\pi}PV\!\!\int_{\R^2}R_1\big( \frac{\dau^4X\cdot N}{|\dau X|^3}\big)(\al)\dau^4X(\al)\cdot\int_{\R^2}\frac{N(\al)\dau\Omega(\al)-  N(\beta)\dau\Omega(\beta)}{|X(\al)-X(\beta)|^3}d\beta d\al,
$$
$$
Q_8=\frac{A_\mu}{4\pi}PV\!\!\int_{\R^2}R_1\big( \frac{\dau^4X\cdot N}{|\dau X|^3}\big)(\al)\int_{\R^2}((\dau\Omega N\cdot\dau^4X)(\al)-(\dau\Omega N\cdot\dau^4X)(\beta))C(\al,\beta)d\beta d\al,
$$
and
$$C(\al,\beta)=\frac{1}{|X(\al)-X(\beta)|^3}-\frac{1}{|\grad X(\al)(\al-\beta)|^3},$$
$$
Q_9=\frac{A_\mu}{4\pi}PV\!\!\int_{\R^2}R_1\big( \frac{\dau^4X\cdot N}{|\dau X|^3}\big)(\al)\frac{1}{|\dau X(\al)|^3}\Lambda(\dau\Omega N\cdot\dau^4X)(\al)d\al.
$$
In $Q_7$ we have
$$
Q_7\leq \|R_1\big(\frac{\dau^4X\cdot N}{|\dau X|^3}\big)\|_{L^2}\|\dau^4X\|_{L^2}\sup_\al \Big|\int_{\R^2}\frac{N(\al)\dau\Omega(\al)-  N(\beta)\dau\Omega(\beta)}{|X(\al)-X(\beta)|^3}d\beta \Big|
$$
giving us the appropriated control, which can be also obtained in $Q_8$ because the corresponding kernel has degree $-2$. Regarding  $Q_9$ we have the expression
\begin{align*}
Q_9&=\frac{A_\mu}{4\pi}PV\!\!\int_{\R^2}R_1\big( \frac{\dau^4X\cdot N}{|\dau X|^3}\big)[\frac{1}{|\dau X|^3}\Lambda(\dau\Omega N\cdot\dau^4X)-\Lambda(\frac{\dau\Omega N\cdot\dau^4X}{|\dau X|^3})]d\al\\
&\quad+\frac{A_\mu}{4\pi}PV\!\!\int_{\R^2}R_1\big( \frac{\dau^4X\cdot N}{|\dau X|^3}\big)\Lambda(\dau\Omega \frac{\dau^4X\cdot N}{|\dau X|^3})d\al.
\end{align*}
Then we   use \eqref{commLambda} to control the first integral above, and since $\Lambda=R_1\dau+R_2\dad$ \eqref{commomega2} we can also take care of the second term.

With $O_5$ one  proceed as we did with $J_6$ \eqref{J6} to get the desired estimate.

Next we use \eqref{oOf} to catch the most singular terms in $O_6$ which are given by
$$
S_7=-\frac{A_\mu}{4\pi}PV\int_{\R^2}\!\!R_1\big(\frac{\dau^4X\cdot N}{|\dau X|^3}\big)(\al)\int_{\R^2}\!\!\frac{(X(\al)-X(\beta))\wedge \dau X(\beta)\cdot\dau X(\al)}{|X(\al)-X(\beta)|^3}\dau^4\dad\Omega(\beta)d\al,
$$
$$
S_8=-\frac{A_\mu}{8\pi^2}PV\int_{\R^2}\!\!R_1\big( \frac{\dau^4X\cdot N}{|\dau X|^3}\big)(\al)\int_{\R^2}\!\!\frac{(X(\al)-X(\beta))\wedge \dau X(\al)}{|X(\al)-X(\beta)|^3}\cdot \dad\Omega(\beta)\dau^5X(\beta)d\al,
$$$$
S_9=\frac{A_\mu}{8\pi^2}PV\int_{\R^2}\!\!R_1\big( \frac{\dau^4X\cdot N}{|\dau X|^3}\big)(\al)\int_{\R^2}\!\!\frac{(X(\al)-X(\beta))\wedge \dad X(\beta)\cdot\dau X(\al)}{|X(\al)-X(\beta)|^3} \dau^5\Omega(\beta)d\al,
$$$$
S_{10}=\frac{A_\mu}{8\pi^2}PV\int_{\R^2}\!\!R_1\big( \frac{\dau^4X\cdot N}{|\dau X|^3}\big)(\al)\int_{\R^2}\!\!\frac{(X(\al)-X(\beta))\wedge \dau X(\al)}{|X(\al)-X(\beta)|^3}\cdot \dau\Omega(\beta)\dau^4\dad X(\beta)d\al.
$$
One may write
$$
S_7=\frac{A_\mu}{4\pi}PV\int_{\R^2}\!\!R_1\big(\frac{\dau^4X\cdot N}{|\dau X|^3}\big)(\al)\int_{\R^2}\!\!\frac{(X(\al)-X(\beta))\wedge (\dau X(\al)-\dau X(\beta))\cdot\dau X(\beta)}{|X(\al)-X(\beta)|^3}\dau^4\dad\Omega(\beta)d\al,
$$
 expressing the fact that we have a kernel of degree $-1$ applied to $\dau^4\dad\Omega$ and, therefore, an integration by parts gives us the desired control as we did before. To treat  $S_8$ we decompose further $S_8=T_1+T_2$:
$$
T_1=-\frac{A_\mu}{4\pi}PV\int_{\R^2}\!\!R_1\big( \frac{\dau^4X\cdot N}{|\dau X|^3}\big)(\al)\int_{\R^2}\!\!D(\al,\beta)\cdot \dad\Omega(\beta)\dau^5X(\beta)d\al,
$$
where $$D(\al,\beta)=\big(\frac{(X(\al)-X(\beta))}{|X(\al)-X(\beta)|^3}-\frac{\grad X(\al)(\al-\beta)}{|\grad X(\al)(\al-\beta)|^3}\big)\wedge \dau X(\al),$$
and
$$
T_2=\frac{A_\mu}{4\pi}PV\int_{\R^2}\!\!R_1\big( \frac{\dau^4X\cdot N}{|\dau X|^3}\big)(\al)\frac{N(\al)}{|\dau X(\al)|^3}\cdot R_2(\dad\Omega\dau^5X)(\al)d\al.
$$
In $T_1$ we use the estimate for the operator \eqref{operadorregular}. The term $T_2$ reads as follows:
\begin{align*}
T_2&=-\frac{A_\mu}{4\pi}PV\int_{\R^2}\!\!R_1\big( \frac{\dau^4X\cdot N}{|\dau X|^3}\big)\frac{N}{|\dau X|^3}\cdot R_2(\dad\dau\Omega\dau^4X)d\al\\
&\quad +\frac{A_\mu}{4\pi}PV\int_{\R^2}\!\!R_1\big( \frac{\dau^4X\cdot N}{|\dau X|^3}\big)[\frac{N}{|\dau X|^3}\cdot (R_2\dau)(\dad\Omega\dau^4X)-(R_2\dau)(\dad\Omega\frac{N\cdot \dau^4X}{|\dau X|^3})]d\al\\
&\quad -\frac{A_\mu}{4\pi}PV\int_{\R^2}\!\!R_1\big( \frac{\dau^4X\cdot N}{|\dau X|^3}\big)(R_2\dau)(\dad\Omega\frac{N\cdot \dau^4X}{|\dau X|^3})d\al.
\end{align*}
The first integral above is easy to estimate, whether for the second one we use \eqref{commderR} and \eqref{commomega2} for the third.

For the next term $S_9$ one has $S_9=T_3+T_4$ where
$$
T_3=\frac{A_\mu}{4\pi}PV\int_{\R^2}\!\!R_1\big( \frac{\dau^4X\cdot N}{|\dau X|^3}\big)(\al)\int_{\R^2}\!\!\frac{(X(\al)-X(\beta))\cdot \dad X(\beta)\wedge(\dau X(\al)-\dau X(\beta))}{|X(\al)-X(\beta)|^3} \dau^5\Omega(\beta)d\al,
$$
$$
T_4=-A_\mu \int_{\R^2}\!\!R_1\big( \frac{\dau^4X\cdot N}{|\dau X|^3}\big)\mathcal{D}(\dau^5\Omega)d\al,
$$
Proceeding as before we get bounds for $T_3$  and the double layer potential estimates help us  to control $T_4$.

 For $S_{10}$ one can adapt exactly  the same approach used for $S_8$. Finally we have to deal with $O_7$ which is given by
$$
O_7=-A_\mu PV\int_{\R^2}BR(X,\omega)\cdot \dau^4 X(R_1\dau)\big( \frac{\dau^4X\cdot N}{|\dau X|^3}\big)d\al,
$$
after an integration by parts. Let us introduce the splitting  $O_7=\sum_{j,k=1}^3 U_j^k$ where
$$
U_j^k=-A_\mu PV\int_{\R^2}BR_j(X,\omega)\dau^4 X_j(R_1\dau)\big(\frac{\dau^4X_k N_k}{|\dau X|^3}\big)d\al.
$$
Then the commutator estimates allows us to write $U_j^k=V_j^k+\mbox{lower order terms}$, where
$$
V_j^k=-A_\mu PV\int_{\R^2}BR_j(X,\omega) \dau^4 X_j\frac{N_k}{|\dau X|^3} (R_1\dau)(\dau^4X_k)d\al.
$$
 Using \eqref{b4a} and \eqref{b4c} one has
$$
N_1\dau^4 X_2=N_2\dau^4 X_1+\mbox{l.o.t.}
$$
so that $V_2^1$ becames
\begin{align*}
\begin{split}
V_2^1=-A_\mu PV\int_{\R^2}\frac{BR_2(X,\omega)N_2}{|\dau X|^3} \dau^4 X_1 (R_1\dau)(\dau^4X_1)d\al-
A_\mu PV\int_{\R^2}f (R_1\dau)(\dau^4X_1)d\al
\end{split}
\end{align*}
where $f$ is at the level of $\dai^3 X$. Integration by parts in the last integral above allows us to conclude that
\begin{equation*}
V_2^1\leq-A_\mu PV\int_{\R^2}\frac{BR_2(X,\omega)N_2}{|\dau X|^3} \dau^4 X_1 (R_1\dau)(\dau^4X_1)d\al +P(|||X|||_{4}).
\end{equation*}
With the help of \eqref{b4a} and \eqref{b4c} we also get
$$
N_1\dau^4 X_3=N_3\dau^4 X_1+\mbox{l.o.t.}
$$
 and therefore
\begin{equation*}
V_3^1\leq-A_\mu PV\int_{\R^2}\frac{BR_3(X,\omega)N_3}{|\dau X|^3} \dau^4 X_1 (R_1\dau)(\dau^4X_1)d\al +P(|||X|||_{4}).
\end{equation*}
Using the two inequalities above we obtain
\begin{equation}\label{V1}
V^1_1+V^1_2+V_3^1\leq-A_\mu PV\int_{\R^2}\frac{BR(X,\omega)\cdot N}{|\dau X|^3} \dau^4 X_1 (R_1\dau)(\dau^4X_1)d\al +P(|||X|||_{4}).
\end{equation}
Next let us observe that
$$
N_2\dau^4X_1=N_1\dau^4X_2+\mbox{l.o.t.},\quad N_2\dau^4X_3=N_3\dau^4X_2+\mbox{l.o.t.},
$$
which implies the estimate
\begin{equation}\label{V2}
V_1^2+V_2^2+V_3^2\leq-A_\mu PV\int_{\R^2}\frac{BR(X,\omega)\cdot N}{|\dau X|^3} \dau^4 X_2 (R_1\dau)(\dau^4X_2)d\al +P(|||X|||_{4}).
\end{equation}
Regarding  $V_1^3$ and $V_2^3$ the identities
$$
N_3\dau^4X_1=N_1\dau^4X_3+\mbox{l.o.t.},\quad N_3\dau^4X_3=N_2\dau^4X_3+\mbox{l.o.t.},
$$
yield
\begin{equation}\label{V3}
V_1^3+V_2^3+V_3^3\leq-A_\mu PV\int_{\R^2}\frac{BR(X,\omega)\cdot N}{|\dau X|^3} \dau^4 X_3 (R_1\dau)(\dau^4X_3)d\al +P(|||X|||_{4}).
\end{equation}
Finally \eqref{V1}, \eqref{V2} and \eqref{V3} imply
$$
\sum_{j,k=1}^3 V_j^k\leq -A_\mu PV\int_{\R^2}\frac{BR(X,\omega)\cdot N}{|\dau X|^3} \dau^4 X\cdot (R_1\dau)(\dau^4X)d\al +P(|||X|||_{4}).
$$
Now we put together all those estimates  (\eqref{M22} - \eqref{V3}) to conclude that
$$
M_2\leq -A_\mu PV\int_{\R^2}\frac{BR(X,\omega)\cdot N}{|\dau X|^3} \dau^4 X\cdot (R_1\dau)(\dau^4X)d\al +P(|||X|||_{4}),
$$
and taking into account  \eqref{densidadRdau} we obtain
\begin{equation}\label{L3}
\widetilde{L}_3=M_1+M_2\leq -\frac{1}{\mu_2\!+\!\mu_1}PV\int_{\R^2}\frac{\sigma}{|\dau X|^3} \dau^4 X\cdot (R_1\dau)(\dau^4X)d\al +P(|||X|||_{4}).
\end{equation}

Finally we have  to work with $L_5$ which can be written in the following manner
$$
L_5=\widetilde{L}_{5}-\frac12PV\int_{\R^2}\dau^4X\cdot \frac{\dad X}{|\dad X|^3}\wedge [R_2(\dau^4\dad\Omega \dau X)-R_2(\dau^4\dad\Omega)\dau X]d\al,
$$
where
$$
\widetilde{L}_5=\frac{1}{2}PV\int_{\R^2}\dau^4X\cdot \frac{N}{|\dad X|^3}(R_2\dad)(\dau^4\Omega)d\al.
$$
Using the commutator estimate, once more, it remains only  to consider $\widetilde{L}_5$, but  let us  point out that replacing the operator $R_1\dau$ by $R_2\dad$ the term $\widetilde{L}_3$ \eqref{Lt3} becomes $\widetilde{L}_5$. Therefore, proceeding exactly as we did  before, one  obtains inequality
\begin{equation}\label{L5}
\widetilde{L}_5\leq -\frac{1}{\mu_2\!+\!\mu_1}PV\int_{\R^2}\frac{\sigma}{|\dau X|^3} \dau^4 X\cdot (R_2\dad)(\dau^4X)d\al +P(|||X|||_{4}).
\end{equation}
Introducing now the identity $\Lambda=(R_1\dau)+(R_2\dad)$  in \eqref{L3} and \eqref{L5} we get
\begin{equation*}
\widetilde{L}_3+\widetilde{L}_5\leq -\frac{1
}{\mu_2\!+\!\mu_1}PV\int_{\R^2}\frac{\sigma}{|\dau X|^3} \dau^4 X\cdot \Lambda(\dau^4X)d\al+P(|||X|||_{4}).
\end{equation*}
Finally all the estimates so far obtained, beginning with  \eqref{dauX}, allow us to write
\begin{equation}\label{final1}
\frac12\frac{d}{dt}\|\dau^4X\|_{L^2}^2(t)\leq-\frac1{\mu_2\!+\!\mu_1}PV\int_{\R^2}\frac{\sigma}{|\dau X|^3} \dau^4 X\cdot \Lambda(\dau^4X)d\al+P(|||X|||_{4}).
\end{equation}
In a similar manner,  using now equations \eqref{fddO},\eqref{b4b} and \eqref{b4d} instead of \eqref{fduO}, \eqref{b4a} and \eqref{b4c} respectively, we obtain
\begin{equation}\label{final2}
\frac12\frac{d}{dt}\|\dad^4X\|_{L^2}^2(t)\leq-\frac1{\mu_2\!+\!\mu_1}PV\int_{\R^2}\frac{\sigma}{|\dau X|^3} \dad^4 X\cdot \Lambda(\dad^4X)d\al+P(|||X|||_{4}).
\end{equation}
Being these two inequalities  \eqref{final1} and \eqref{final2}  the main purpose of this section.


\section{Estimates for the evolution of $\|F(X)\|_{L^\infty}$ and R-T.}


In this section we analyze the evolution of the non-selfintersecting condition of the free surface as well as the Rayleigh-Taylor property, but in order to do that we shall need precise bounds for both $\grad X_t$ and $\Omega_t$.

We shall  estimate $\|\grad X_t\|_{H^k}$ by means of equality \eqref{eqevcon} to get
\begin{equation}\label{qay2}
\|\grad X_t\|_{H^k}\leq P(\|X\|^2_{k+2}+\|F(X)\|^2_{L^\infty}+\||N|^{-1}\|_{L^\infty}),
\end{equation}
for $k\geq 2$. In fact
$$
\|\grad X_t\|_{H^k}\leq \|\grad BR(X,\omega)\|_{H^k}+\|\grad(C_1\dau X+C_2\dad X)\|_{H^k}
$$
and with the help of \eqref{estimateBR} we can handle both terms on the right.

Next we shall consider the norms $\|\Omega_t\|_{H^k}$ to obtain the inequality
\begin{equation}\label{qay3}
\|\Omega_t\|_{H^k}\leq P(\|X\|^2_{k+1}+\|F(X)\|^2_{L^\infty}+\||N|^{-1}\|_{L^\infty}),
\end{equation}
for $k\geq 3$. 
To do that let us  take a time derivative in the identity \eqref{eqOmega} to get
\begin{equation*}\label{ehlp}
\Omega_t(\al,t)-A_\mu\mathcal{D}(\Omega_t)(\al,t)=A_\mu I_1(\al,t)-2A_\rho \partial_tX_{3}(\al,t),
\end{equation*}
which yields
$$
\|\Omega_t\|_{H^1}\leq C\|(I-A_\mu\mathcal{D})^{-1}\|_{H^1}(\|I_1\|_{H^1}+\|\partial_tX_{3}\|_{H^1}),
$$
and since we have control of $\|(I-A_\mu\mathcal{D})^{-1}\|_{H^1}$ and $\|\partial_tX_{3}\|_{H^1}$ it only remains to estimate $\|I_1\|_{H^1}$. For that purpose let us  consider the splitting $I_1=J_1+J_2+J_3$ where
$$
J_1=\frac{1}{2\pi}PV\int_{\R^2}\frac{X_t(\al)-X_t(\al-\beta)}{|X(\al)-X(\al-\beta)|^3}\cdot N(\al-\beta) \Omega(\al-\beta)d\beta,
$$
$$
J_2=\frac{-3}{4\pi}\int_{\R^2}(X(\al)-X(\al-\beta))\cdot(X_t(\al)-X_t(\al-\beta))\frac{X(\al)-X(\al-\beta)}{|X(\al)-X(\al-\beta)|^5}\cdot N(\al-\beta) \Omega(\al-\beta)d\beta,$$
$$
J_3=\frac{1}{2\pi}PV\int_{\R^2}\frac{X(\al)-X(\al-\beta)}{|X(\al)-X(\al-\beta)|^3}\cdot N_t(\al-\beta) \Omega(\al-\beta)d\beta.
$$
Proceeding as we did with the  operator $\mathcal{T}_2$ \eqref{operadorgeneral2} ( with $X_t$ instead of $\daj X_k$) one get
$$
\|J_1\|_{L^2}+\|J_2\|_{L^2}\leq P(\|X\|_{4}+\|F(X)\|_{L^\infty}+\||N|^{-1}\|_{L^\infty}).
$$
Regarding $J_3$ we split further
$$
J_3=\frac{1}{2\pi}\int_{|\beta|>1}d\beta+\frac{1}{2\pi}\int_{|\beta|<1}d\beta=K_1+K_2.
$$
Since
$$
|K_1(\al)|\leq \|F(X)\|^2_{L^\infty}\int_{|\beta|>1}\frac{|N_t(\al-\beta)||\Omega(\al-\beta)|}{2\pi|\beta|^2}d\beta,
$$
Young's inequality yields
$$
\|K_1\|_{L^2}\leq \|F(X)\|^2_{L^\infty}\|N_t\Omega\|_{L^1}\leq C\|F(X)\|^2_{L^\infty}\|N_t\|_{L^2}\|\Omega\|_{L^{2}},
$$
and since we know that $\|N_t\|_{L^2}\leq \|\grad X\|_{L^\infty}\|\grad X_t\|_{L^2}$,  estimate \eqref{qay2} allows us to handle the terms  $K_1$.
The estimate for $K_2$ is similar to the one obtained for $I_2$ \eqref{dT1} in the appendix.

 Next we consider the most singular terms in $\dau I_1$ which are given by
$$
J_4=\frac{1}{2\pi}PV\int_{\R^2}\frac{\dau X_t(\al)-\dau X_t(\al-\beta)}{|X(\al)-X(\al-\beta)|^3}\cdot N(\al-\beta) \Omega(\al-\beta)d\beta,
$$
$$
J_5=\frac{-3}{4\pi}\int_{\R^2}(X(\al)-X(\al-\beta))\cdot(\dau X_t(\al)-\dau X_t(\al-\beta))\frac{X(\al)-X(\al-\beta)}{|X(\al)-X(\al-\beta)|^5}\cdot N(\al-\beta) \Omega(\al-\beta)d\beta,$$
$$
J_6=\frac{1}{2\pi}PV\int_{\R^2}\frac{X(\al)-X(\al-\beta)}{|X(\al)-X(\al-\beta)|^3}\cdot \dau N_t(\al-\beta) \Omega(\al-\beta)d\beta.
$$
because the remainder terms are easier to handle. Let us write $J_4=K_3+K_4$ where
$$
K_3=\frac{1}{2\pi}PV\int_{\R^2}\frac{\dau X_t(\al)-\dau X_t(\al-\beta)}{|X(\al)-X(\al-\beta)|^3}\cdot (N(\al-\beta) \Omega(\al-\beta)-N(\al) \Omega(\al))d\beta,
$$$$
K_4=\frac{1}{2\pi}PV\int_{\R^2}\frac{\dau X_t(\al)-\dau X_t(\al-\beta)}{|X(\al)-X(\al-\beta)|^3}\cdot N(\al) \Omega(\al)d\beta.
$$
In $K_3$, the identity $\dau X_t(\al)-\dau X_t(\al-\beta)=\int_0^1\grad\dau X_t(\al+(s-1)\beta)ds\cdot \beta$ together with \eqref{qay2} gives us the desired control. Regarding  $K_4$ we may observe its similarity with $\mathcal{T}_3$ \eqref{operadorgeneral3}, so that an application to \eqref{qay2} yields the appropriated bound;  $J_5$ can be treated in a similar manner and $J_6$ is analogous to  $J_3$. By symmetry, one could get the same estimate for $\dad I_1$, so that  finally:
\begin{equation}\label{qay4}
\|\Omega_t\|_{H^1}\leq P(\|X\|^2_{4}+\|F(X)\|^2_{L^\infty}+\||N|^{-1}\|_{L^\infty}).
\end{equation}

Next, we will show how to deal with $\|\Omega_t\|_{H^2}$. Using equation \eqref{fduO} one gets
$$
\dau^2\Omega_t=-2A_\mu
\dau\dpt (BR(X,\omega)\cdot\dau X)-2A_\rho\dau^2\dpt X_3,
$$
and with the help of \eqref{qay2}, the last term above is properly controlled.
To continue we shall consider the most singular remainder  terms. Namely, in
$-\dau\dpt (BR(X,\omega)\cdot\dau X)$, we have:
$$
L_1=-BR(X,\omega)\cdot\dau^2 X_t
$$
$$L_2=\frac{1}{4\pi}PV\int_{\R^2}\frac{\dau X_t(\al)-\dau X_t(\al-\beta)}{|X(\al)-X(\al-\beta)|^3}\wedge \omega(\al-\beta)d\beta
\cdot \dau X(\al)$$
$$L_3=\frac{-3}{8\pi}PV\int_{\R^2}A(\al,\beta)\frac{X(\al)-X(\al-\beta)}{|X(\al)-X(\al-\beta)|^5}\wedge \omega(\al-\beta)d\beta
\cdot \dau X(\al)$$
where $A(\al,\beta)=(X(\al)-X(\al-\beta))\cdot(\dau X_t(\al)-\dau X_t(\al-\beta))$,
$$L_4=\frac{1}{2\pi}PV\int_{\R^2}\frac{X(\al)-X(\al-\beta)}{|X(\al)-X(\al-\beta)|^3}\wedge \dau\omega_t(\al-\beta)d\beta
\cdot \dau X(\al).$$
 Let us observe that $\|L_1\|_{L^2}\leq \|BR(X,\omega)\|_{L^\infty}\|\dau^2 X_t\|_{L^2}$, where both  quantities have been appropriately controlled  before. In $L_2$ and $L_3$ we have kernels of degree  $-2$, and therefore operators analogous to $\mathcal{T}_3$ \eqref{operadorgeneral3} acting on $\dau X_t$. Therefore  using \eqref{qay2} its control follows easily. In $L_4$ we use the decomposition
$$
L_4=\frac{1}{2\pi}PV\int_{|\beta|>1}d\beta+\frac{1}{2\pi}PV\int_{|\beta|<1}d\beta=M_1+M_2.
$$
Thus an integration by parts yields
$$
\|M_1\|_{L^2}\leq C\|F(X)\|^3_{L^\infty}\|\grad X\|^2_{L^\infty}\|w_t\|_{L^2}.
$$
Formula \eqref{oOf} together with estimates \eqref{qay2} and \eqref{qay4}, provides the appropriated bound.

Next let us  expand \eqref{oOf} to obtain  the most singular terms in $M_2$ which are given by the integrals:
$$
O_1=-\frac{A_\mu}{2\pi}PV\int_{|\beta|<1}
\frac{X(\al)-X(\al-\beta)}{|X(\al)-X(\al-\beta)|^3}\wedge\dad\Omega(\al-\beta)\dau^2 X_t(\al-\beta)d\beta\cdot\dau X(\al),
$$
$$
O_2=-\frac{A_\mu}{2\pi}PV\int_{|\beta|<1}
\frac{X(\al)-X(\al-\beta)}{|X(\al)-X(\al-\beta)|^3}\wedge\dau\dad\Omega_t(\al-\beta)\dau X(\al-\beta)d\beta\cdot\dau X(\al),
$$$$
O_3=\frac{A_\mu}{2\pi}PV\int_{|\beta|<1}
\frac{X(\al)-X(\al-\beta)}{|X(\al)-X(\al-\beta)|^3}\wedge\dau\Omega(\al-\beta)\dau\dad X_t(\al-\beta)d\beta\cdot\dau X(\al),
$$$$
O_4=\frac{A_\mu}{2\pi}PV\int_{|\beta|<1}
\frac{X(\al)-X(\al-\beta)}{|X(\al)-X(\al-\beta)|^3}\wedge\dau^2\Omega_t(\al-\beta)\dad X(\al-\beta)d\beta\cdot\dau X(\al).
$$
Estimate \eqref{qay2} help us with the terms  $O_1$ and $O_3$, which can be treated with  the same approach used for $I_2$ \eqref{dT1} in the appendix. Let us  write $O_2$ as follows
$$
O_2=\frac{A_\mu}{2\pi}\!\!\int_{|\beta|<1}
\frac{X(\al)\!-\!X(\al\!-\!\beta)}{|X(\al)\!-\!X(\al\!-\!\beta)|^3}\wedge\dau\dad\Omega_t(\al\!-\!\beta)(\dau X(\al)\!-\!\dau X(\al\!-\!\beta))d\beta\cdot\dau X(\al),
$$
which can be estimated integrating by parts in the variable $\beta_1$ using the following identity
$$\dau\dad\Omega_t(\al\!-\!\beta)=-\dbu(\dad\Omega_t(\al\!-\!\beta)).$$ Let us point
out that the kernel in the integral $O_2$ has degree $-1$ and, therefore, one can use  \eqref{qay4} to control it. It remains to deal with $O_4$ which is decomposed in the form  $O_4=P_1+P_2$, where
$$
P_1=\frac{A_\mu}{2\pi}PV\!\!\int_{|\beta|<1}
\frac{X(\al)\!-\!X(\al\!-\!\beta)}{|X(\al)\!-\!X(\al\!-\!\beta)|^3}\wedge\dau^2\Omega_t(\al\!-\!\beta)(\dad X(\al\!-\!\beta)\!-\!\dad X(\al))d\beta\cdot\dau X(\al),
$$
$$
P_2=-\frac{A_\mu}{2\pi}PV\int_{|\beta|<1}
\frac{X(\al)-X(\al-\beta)}{|X(\al)-X(\al-\beta)|^3}\dau^2\Omega_t(\al-\beta)d\beta\cdot N(\al).
$$
 $P_1$ is estimated like  $O_2$. We rewrite $P_2$ as follows
$$
P_2=-\frac{A_\mu}{2\pi}PV\int_{|\beta|<1}\Big(
\frac{X(\al)-X(\al-\beta)}{|X(\al)-X(\al-\beta)|^3}-\frac{\grad X(\al)\cdot\beta}{|\grad X(\al)\cdot\beta|^3}\Big)\dau^2\Omega(\al-\beta)d\beta\cdot N(\al),
$$
 and this expression shows that the above integral can be estimated like $\mathcal{T}_4$ \eqref{operadorregular}.

  Using \eqref{qay4} we obtain
\begin{equation*}
\|\dau^2\Omega_t\|_{L^2}\leq P(\|X\|^2_{4}+\|F(X)\|^2_{L^\infty}+\||N|^{-1}\|_{L^\infty}),
\end{equation*}
and the identity
$$
\dad^2\Omega_t=-2A_\mu
\dad\dpt (BR(X,\omega)\cdot\dad X)-2A_\rho\dad^2\dpt X_3,
$$
yields
\begin{equation*}\label{nose}
\|\dad^2\Omega_t\|_{L^2}\leq P(\|X\|^2_{4}+\|F(X)\|^2_{L^\infty}+\||N|^{-1}\|_{L^\infty}),
\end{equation*}
that is:
\begin{equation}\label{qay5}
\|\Omega_t\|_{H^2}\leq P(\|X\|^2_{4}+\|F(X)\|^2_{L^\infty}+\||N|^{-1}\|_{L^\infty}).
\end{equation}
Next we consider third order derivatives
$$
\dau^3\Omega_t=-2A_\mu
\dau^2\dpt (BR(X,\omega)\cdot\dau X)-2A_\rho\dau^3\dpt X_3.
$$
Since \eqref{qay2} gives us control of the last term, we will concentrate in the other one which is a much more diffecult character. In particular, for
$-\dau^2\dpt (BR(X,\omega)\cdot\dau X)$, the most singular component are given by
$$
L_5=-BR(X,\omega)\cdot\dau^3 X_t
$$
$$L_6=\frac{1}{4\pi}PV\int_{\R^2}\frac{\dau^2 X_t(\al)-\dau^2 X_t(\al-\beta)}{|X(\al)-X(\al-\beta)|^3}\wedge \omega(\al-\beta)d\beta
\cdot \dau X(\al)$$
$$L_7=\frac{-3}{8\pi}PV\int_{\R^2}B(\al,\beta)\frac{X(\al)-X(\al-\beta)}{|X(\al)-X(\al-\beta)|^5}\wedge \omega(\al-\beta)d\beta
\cdot \dau X(\al)$$
where $B(\al,\beta)=(X(\al)-X(\al-\beta))\cdot(\dau^2 X_t(\al)-\dau^2 X_t(\al-\beta))$,
$$L_8=\frac{1}{2\pi}PV\int_{\R^2}\frac{X(\al)-X(\al-\beta)}{|X(\al)-X(\al-\beta)|^3}\wedge \dau^2\omega_t(\al-\beta)d\beta
\cdot \dau X(\al).$$
Inequalities \eqref{qay2} and \eqref{qay5} show how to handle $L_i$, $i=5,...,8$ as $L_j$, $j=1,...,4$ respectively, then a  similar approach for $\dad^3\Omega_t$ allows us to get finally \eqref{qay3} for k=3. The cases $k>3$ are similar to deal with.
\newline

Our next goal is to obtain  estimates for the evolution of $\|F(X)\|_{L^\infty}$ and R-T.
Regarding the quantity $F(X)$ we have
\begin{align}
\begin{split}\label{mct}
\D\dt F(X)(\al,\beta,t)&=-\frac{|\beta|(X(\al,t)-X(\al-\beta,t))\cdot(X_t(\al,t)-
X_t(\al-\beta,t))}{|X(\al,t)-X(\al-\beta,t)|^3}\\
&\leq (F(X)(\al,\beta,t))^2\|\grad X_t\|_{L^\infty}(t).
\end{split}
\end{align}
Then Sobolev inequalities in $\|\grad X_t\|_{L^\infty}(t)$ together with \eqref{qay2} yield
\begin{align*}
\D\dt F(X)(\al,\beta,t)&\leq F(X)(\al,\beta,t)
P(\|X\|^2_{4}(t)+\|F(X)\|^2_{L^\infty}(t)+\||N|^{-1}\|_{L^\infty}(t)),
\end{align*}
and an integration in time gives us
\begin{align*}
F(X)(\al,\beta,t+h)&\leq F(X)(\al,\beta,t)\exp \Big(\int_t^{t+h}P(s)
ds\Big),
\end{align*}
for $h>0$, where $$P(s)=P(\|X\|^2_{4}(s)+\|F(X)\|^2_{L^\infty}(s)+\||N|^{-1}\|_{L^\infty}(s)).$$ Hence
\begin{align*}
\|F(X)\|_{L^\infty}(t+h)&\leq \|F(X)\|_{L^\infty}(t)\exp \Big(\int_t^{t+h}P(s)ds\Big).
\end{align*}
The inequality above applied to the  limit:
$$
\dt \|F(X)\|_{L^\infty}(t)=\lim_{h\rightarrow0^+}\frac{\|F(X)\|_{L^\infty}(t+h)-\|F(X)\|_{L^\infty}(t)}{h}
$$
allows us to get
\begin{align*}
\begin{split}\label{mct2}
\D\dt \|F(X)\|_{L^\infty}(t)\leq \|F(X)\|_{L^\infty}(t)
P(\|X\|^2_{4}+\|F(X)\|^2_{L^\infty}+\||N|^{-1}\|_{L^\infty}).
\end{split}
\end{align*}

Next we search for an a priori estimate for the evolution of
the infimum of the difference of the gradients of the pressure in
the normal direction to the interface. Let us recall the formula
\begin{align*}
\begin{split}
\sigma(\al,t)&=(\mu^2-\mu^1)BR(X,\omega)(\al,t)\cdot N(\al,t)+(\rho^2-\rho^1)N_3(\al,t),
\end{split}
\end{align*}
to obtain
$$\D\dt(\frac{1}{\sigma(\al,t)})=-\frac{\sigma_t(\al,t)}{\sigma^2(\al,t)}$$
with $\sigma_t(\al,t)=I_1+I_2$ where
$$
I_1=((\mu^2-\mu^1)BR(X,\omega)(\al,t)+(\rho^2-\rho^1)(0,0,1))\cdot N_t(\al,t),
$$
$$
I_2=(\mu^2-\mu^1)BR_t(X,\omega)(\al,t)\cdot N(\al,t).
$$
First we  deal with $\|I_1\|_{L^\infty}$ using the estimates \eqref{qay2} for $\grad X_t$,  and then we focus our attention on $I_2$ using the splitting $I_2=J_1+J_2+J_3$ where
$$
J_1=-\frac{1}{4\pi}PV\int_{\R^2}\frac{X_t(\al)-X_t(\al-\beta)}{|X(\al)-X(\al-\beta)|^3}\wedge \omega(\al-\beta)d\beta,
$$
$$
J_2=\frac{3}{4\pi}PV\!\!\int_{\R^2}\!(X(\al)\!-\!X(\al\!-\!\beta))\!\wedge\!\omega(\al\!-\!\beta)\frac{(X(\al)\!-\!X(\al\!-\!\beta))\!\cdot\!(X_t(\al)\!-\!X_t(\al\!-\!\beta))}{|X(\al)\!-\!X(\al\!-\!\beta)|^5}d\beta
$$
$$
J_3=-\frac{1}{4\pi}PV\int_{\R^2}\frac{X(\al)-X(\al-\beta)}{|X(\al)-X(\al-\beta)|^3}\wedge \omega_t(\al-\beta) d\beta.
$$
The terms $J_1$ and  $J_2$ are similar and can be treated with the same method. Let us consider $J_1=K_1+K_2+K_3+K_4$ where
$$
K_1=-\frac{1}{4\pi}\int_{|\beta|>1}\frac{X_t(\al)-X_t(\al-\beta)}{|X(\al)-X(\al-\beta)|^3}\wedge \omega(\al-\beta)d\beta,
$$                                                  $$
K_2=\frac{1}{4\pi}\int_{|\beta|<1}\frac{X_t(\al)-X_t(\al-\beta)}{|X(\al)-X(\al-\beta)|^3}\wedge (\omega(\al)-\omega(\al-\beta))d\beta,
$$$$
K_3=-\frac{1}{4\pi}\int_{|\beta|<1}[\frac{1}{|X(\al)\!-\!X(\al\!-\!\beta)|^3}\!-\!\frac{1}{|\grad X(\al)\!\cdot\!\beta|^3}](X_t(\al)\!-\!X_t(\al\!-\!\beta))\wedge \omega(\al)d\beta,
$$                                                  $$
K_4=-\frac{1}{4\pi}PV\int_{|\beta|<1}\frac{X_t(\al)-X_t(\al-\beta)}{|\grad X(\al)\cdot\beta|^3}\wedge \omega(\al)d\beta,
$$                                                  First we have
$$
\|K_1\|_{L^\infty}\leq C\|F(X)\|^3_{L^\infty}\|\grad X_t\|_{L^\infty}\|\omega\|_{L^2}\big(\int_{|\beta|>1}|\beta|^{-4}d\beta\big)^{1/2}
$$
giving us an appropriated  control. Next, we get
$$
\|K_2\|_{L^\infty}\leq C\|F(X)\|^3_{L^\infty}\|\grad X_t\|_{L^\infty}\|\grad\omega\|_{L^\infty}\int_{|\beta|<1}|\beta|^{-1}d\beta,
$$ and an analogous estimate  for $K_3$. Therefore, Sobolev's embedding  help us to obtain the desired control. Regarding  $K_4$ we have
$$
K_4=-\frac{1}{4\pi}\int_{|\beta|<1}\frac{X_t(\al)-X_t(\al-\beta)-\grad X_t(\al)\cdot\beta}{|\grad X(\al)\cdot\beta|^3}\wedge \omega(\al)d\beta.
$$
Inequality \eqref{sioe} yields
$$
\|K_4\|_{L^\infty}\leq C\|\grad X\|^3_{L^\infty}\||N|^{-1}\|^3_{L^\infty}\|\omega\|_{L^\infty}\|\grad X_t\|_{C^\delta}\int_{|\beta|<1}|\beta|^{-2+\delta}d\beta,
$$
and the control $\|\grad X_t\|_{C^\delta}$ follows again by \eqref{qay2} and
Sobolev's embedding. Next let us continue with $J_3=K_5+K_6$ where
$$
K_5=\frac{-1}{4\pi}PV\!\!\int_{|\beta|>1}\frac{X(\al)\!-\!X(\al\!-\!\beta)}{|X(\al)\!-\!X(\al\!-\!\beta)|^3}\wedge (\dbu((\Omega\dad X)_t(\al\!-\!\beta))\!-\!\dbd((\Omega\dau X)_t(\al\!-\!\beta)))d\beta,
$$
$$
K_6=-\frac{1}{4\pi}PV\int_{|\beta|<1}\frac{X(\al)-X(\al-\beta)}{|X(\al)-X(\al-\beta)|^3}\wedge \omega_t(\al-\beta) d\beta,
$$
Integration by parts yields
$$
\|K_5\|_{L^\infty}\leq C\|F(X)\|^3_{L^\infty}\|\grad X\|_{L^\infty}(\|\Omega\|_{L^\infty}\|\grad X_t\|_{L^\infty}+\|\Omega_t\|_{L^\infty}\|\grad X\|_{L^\infty}),
$$
where $4\pi C=\int_{|\beta|>1}|\beta|^{-3}d\beta+\int_{|\beta|=1}dl(\beta)$, and we may use \eqref{qay3} to estimate $\|\Omega_t\|_{L^\infty}$. With $K_6$ we introduce  a similar splitting to obtain
$$
\|K_6\|_{L^\infty}\leq P(\|X-(\al,0)\|_{C^2}+\|F(X)\|_{L^\infty}+\||N|^{-1}\|_{L^\infty})\|\omega_t\|_{C^\delta}.
$$
 Then it remains to estimate $\|\omega_t\|_{C^\delta}$, for which purpose  we use formula \eqref{oOf} and inequalities \eqref{qay2}\eqref{qay3}. Therefore we have the estimate: $$\D\dt(\frac{1}{\sigma(\al,t)})\leq \frac{1}{\sigma^2(\al,t)}P(\|X\|_{4}(t)+\|F(X)\|_{L^\infty}(t)+\||N|^{-1}\|_{L^\infty}(t)),$$
and proceeding similarly as we did for $F(X)$ we get finally:
$$\D\dt\|\sigma^{-1}\|_{L^\infty}(t)\leq \|\sigma^{-1}\|^2_{L^\infty}(t)P(\|X\|_{4}(t)+\|F(X)\|_{L^\infty}(t)+\||N|^{-1}\|_{L^\infty}(t)).$$




\section{Appendix}


Here we prove first some helpful inequalities regarding commutators of  Riesz transform  ($R_j$, $j=1,2$) with several differential operators. Next we analyze the  singular integral operators associated to the non-selfintersecting surface which appears throughout the paper. But the main goal of this section, however,  is to simplify the presentation of the main result.

\begin{lemma}
Consider $f\in L^2(\R^2)$, and $g\in C^{1,\delta}(\R^2)$ with $0<\delta<1$. Then for any $k,l=1,2$ we have the following estimate
\begin{equation}\label{commderR}
\|(R_k\partial_{\alpha_l})(gf)-g(R_k\partial_{\alpha_l})(f)\|_{L^2}\leq C\|g\|_{C^{1,\delta}}\|f\|_{L^2}.
\end{equation}
\end{lemma}
An application of the above inequalities to the operator $\Lambda=(R_1\dau)+(R_2\dad)$ yields
\begin{equation}\label{commLambda}
\|\Lambda(gf)-g\Lambda(f)\|_{L^2}\leq C\|g\|_{C^{1,\delta}}\|f\|_{L^2}.
\end{equation}
For vector fields we have
\begin{lemma}
Consider $f,g:\R^2\rightarrow\R^3$ vector fields where $f\in L^2(\R^2)$ and $g\in C^{1,\delta}(\R^2)$ with $0<\delta<1$. Then for any $k,l=1,2$ the following inequality holds
\begin{equation}\label{commomega}
\big|\int_{\R^2}(g\wedge f)\cdot(R_k\partial_{\alpha_l})(f)d\alpha\big|\leq C\|g\|_{C^{1,\delta}}\|f\|^2_{L^2}.
\end{equation}
\end{lemma}

Proof: Denoting  with $I$ the integral above and since the operator $R_k\partial_{\alpha_l}$ is self-adjoint we may write
\begin{align*}
I=&\int_{\R^2}f_1[(R_k\partial_{\alpha_l})(g_2f_3)-g_2(R_k\partial_{\alpha_l})(f_3)]d\alpha+\int_{\R^2}f_2[(R_k\partial_{\alpha_l})(g_3f_1)-g_3(R_k\partial_{\alpha_l})(f_1)]d\alpha\\
&+\int_{\R^2}f_3[(R_k\partial_{\alpha_l})(g_1f_2)-g_1(R_k\partial_{\alpha_l})(f_2)]d\alpha.
\end{align*}
Then estimate \eqref{commderR} yields \eqref{commomega}.

\begin{lemma}
Consider $f\in L^2(\R^2)$ and $g\in C^{1,\delta}(\R^2)$ with $0<\delta<1$. Then for any $j,k,l=1,2$ the following inequality holds
\begin{equation}\label{commomega2}
\big|\int_{\R^2}R_j(f)(R_k\partial_{\alpha_l})(gf)d\alpha\big|\leq C\|g\|_{C^{1,\delta}}\|f\|^2_{L^2}.
\end{equation}
\end{lemma}

Proof: Let $J$ be the integral to be bounded, then we have
\begin{align*}
J&=\int_{\R^2}R_j(f)[(R_k\partial_{\alpha_l})(gf)-g(R_k\partial_{\alpha_l})(f)]d\alpha-
\int_{\R^2}[R_j(fg)-gR_j(f)](R_k\partial_{\alpha_l})(f)d\alpha\\
&\quad +\int_{\R^2}R_j(fg)(R_k\partial_{\alpha_l})(f)d\alpha
\end{align*}
Since $R^{\ast}_j = -R_j$  and $R_k\partial_{\alpha_l}$ is self-adjoint we get
\begin{align*}
J&=\frac12\int_{\R^2}R_j(f)[(R_k\partial_{\alpha_l})(gf)-g(R_k\partial_{\alpha_l})(f)]d\alpha-
\frac12\int_{\R^2}[R_j(fg)-gR_j(f)](R_k\partial_{\alpha_l})(f)d\alpha.
\end{align*}
An integration by parts in the second integral above yields
\begin{align*}
J&=\frac12\int_{\R^2}R_j(f)[(R_k\partial_{\alpha_l})(gf)-g(R_k\partial_{\alpha_l})(f)]d\alpha+
\frac12\int_{\R^2}[(R_j\partial_{\alpha_l})(fg)-g(R_j\partial_{\alpha_l})(f)](R_k)(f)d\alpha\\
&\quad -\frac12\int_{\R^2}(\partial_{\alpha_l}g)R_j(f)R_k(f)d\alpha,
\end{align*}
allowing  us to conclude the proof.

\begin{lemma}
 Let us  define for any $j=1,2$ and $k=1,2,3$ the following operators:
\begin{equation}\label{operadorgeneral}
\mathcal{T}_1(\partial_{\alpha_j}f)(\al)=PV\int_{\R^2}\frac{X_k(\al)-X_k(\al-\beta)}{|X(\al)-X(\al-\beta)|^3}\partial_{\alpha_j}f(\al-\beta)d\beta,
\end{equation}
\begin{equation}\label{operadorgeneral2}
\mathcal{T}_2(f)(\al)=PV\int_{\R^2}\frac{\partial_{\alpha_j}X_k(\al)-\partial_{\alpha_j}X_k(\al-\beta)}{|X(\al)-X(\al-\beta)|^3}f(\al-\beta)d\beta,
\end{equation}
\begin{equation}\label{operadorgeneral3}
\mathcal{T}_3(f)(\al)=PV\int_{\R^2}\frac{f(\al)-f(\al-\beta)}{|X(\al)-X(\al-\beta)|^3}d\beta,
\end{equation}
\begin{equation}\label{operadorregular}
\mathcal{T}_4(\partial_{\alpha_j} f)(\al)=PV\int_{\R^2}
\Big(\frac{(X(\al)-X(\beta))}{|X(\al)-X(\beta)|^3}-\frac{\grad X(\al)\cdot(\al-\beta)}{|\grad X(\al)\cdot(\al-\beta)|^3}\Big) \partial_{\alpha_j} f(\beta)d\beta d\al,
\end{equation}
where $\grad X(\al)\cdot\beta=\partial_{\alpha_1}X(\al)\beta_1+\partial_{\alpha_2} X(\al)\beta_2$.
Assume that $X(\al)-(\al,0)\in C^{2,\delta}(\R^2)$, and that both $F(X)$ and $|N|^{-1}$ are in $ L^\infty$ where $$F(X)(\al,\beta)=|\beta|/|X(\al)-X(\al-\beta)|\quad \mbox{and}\quad N(\al)=\dau X(\al) \wedge \dad X(\al).$$ Then the following estimates hold:
\begin{equation}\label{estimaciongeneral}
\|\mathcal{T}_1(\partial_{\alpha_j}f)\|_{L^2}\leq P(\|X-(\al,0)\|_{C^{1,\delta}}+\|F(X)\|_{L^\infty}+\||N|^{-1}\|_{L^\infty})(\|f\|_{L^2}+\|\partial_{\alpha_j}f\|_{L^2}),
\end{equation}
\begin{equation}\label{estimaciongeneral2}
\|\mathcal{T}_2(f)\|_{L^2}\leq P(\|X-(\al,0)\|_{C^{2,\delta}}+\|F(X)\|_{L^\infty}+\||N|^{-1}\|_{L^\infty})\|f\|_{L^2},
\end{equation}
\begin{equation}\label{estimaciongeneral3}
\|\mathcal{T}_3(f)\|_{L^2}\leq P(\|X-(\al,0)\|_{C^{2,\delta}}+\|F(X)\|_{L^\infty}+\||N|^{-1}\|_{L^\infty})\|f\|_{H^1},
\end{equation}
\begin{equation}\label{estimacionregular}
\|\mathcal{T}_4(f)\|_{L^2}\leq P(\|X-(\al,0)\|_{C^{2,\delta}}+\|F(X)\|_{L^\infty}+\||N|^{-1}\|_{L^\infty})\|f\|_{L^2},
\end{equation}
with $P$ a polynomial function.
\end{lemma}

Proof: To estimate the first set of operators we consider first the splitting
\begin{equation}\label{dT1}
\mathcal{T}_1(\partial_{\alpha_j}f)=PV\int_{|\beta|>1}d\beta+PV\int_{|\beta|<1}d\beta=I_1+I_2
\end{equation}
and an integration by parts allows us to write $I_1=J_1+J_2+J_3$ where
$$
J_1=\int_{|\beta|>1}\frac{-\partial_{\alpha_j}X_k(\al-\beta)}{|X(\al)-X(\al-\beta)|^3}f(\al-\beta)d\beta,
$$
$$
J_2=3\int_{|\beta|>1}\frac{(X_k(\al)-X_k(\al-\beta))(X(\al)-X(\al-\beta))\cdot \partial_{\alpha_j}
X(\al-\beta)}{|X(\al)-X(\al-\beta)|^5}f(\al-\beta)d\beta,
$$
and
$$
J_3=\int_{|\beta|=1}\frac{X_k(\al)-X_k(\al-\beta)}{|X(\al)-X(\al-\beta)|^3}f(\al-\beta)dl(\beta).
$$
The above decomposition shows that
$$
|I_1|\leq C\|X-(\al,0)\|_{C^{1}}\|F(X)\|^3_{L^\infty}(\int_{|\beta|>1}\frac{|f(\al-\beta)|}{|\beta|^3}d\beta+\int_{|\beta|=1}|f(\al-\beta)|dl(\beta))
$$ and then Minkowski's inequality gives the desired control.

Regarding  $I_2$ we write $I_2=J_4+J_5+J_6$ with
$$
J_4=\int_{|\beta|<1}\frac{X_k(\al)-X_k(\al-\beta)-\grad X_k(\al)\cdot\beta}{|X(\al)-X(\al-\beta)|^3}\partial_{\alpha_j}f(\al-\beta)d\beta,
$$
$$
J_5=\grad X_k(\al)\cdot\int_{|\beta|<1}\beta [\frac{1}{|X(\al)-X(\al-\beta)|^3}-\frac{1}{|\grad X(\al)\cdot \beta|^3}]\partial_{\alpha_j}f(\al-\beta)d\beta,
$$
$$
J_6=\grad X_k(\al)\cdot PV\int_{|\beta|<1}\frac{\beta}{|\grad X(\al)\cdot\beta|^3}\partial_{\alpha_j}f(\al-\beta)d\beta.
$$

It is easy to see that
\begin{equation}\label{ABJ4}
J_4\leq \|X-(\al,0)\|_{C^{1,\delta}}\|F(X)\|^3_{L^\infty}\int_{|\beta|<1}\frac{|\partial_{\alpha_j}f(\al-\beta)|}{|\beta|^{2-\delta}}d\beta,
\end{equation}
and therefore that term can be estimated also with the use of  Minkowski's inequality.

Some elementary algebraic manipulations allows us to get
$$
J_5\leq C\|X-(\al,0)\|^2_{C^{1,\delta}}\int_{|\beta|<1}[(F(X)(\al,\beta))^4+\frac{|\beta|^4}{|\grad X(\al)\cdot \beta|^4}]\frac{|\partial_{\alpha_j}f(\al-\beta)|}{|\beta|^{2-\delta}}d\beta,
$$
and then the inequality
\begin{equation}\label{sioe}\frac{|\beta|}{|\grad X(\al)\cdot \beta|}\leq 2 \|\grad X\|_{L^\infty}\||N|^{-1}\|_{L^\infty}\end{equation}
yields for $J_5$ the same estimate \eqref{ABJ4}.

The term $J_6$ can be written as
$$
J_6=\grad X_k(\al)\cdot PV\int_{|\beta|<1}\frac{\Sigma(\al,\beta)}{|\beta|^2}\partial_{\alpha_j}f(\al-\beta)d\beta,
$$
where
$$
(i)\,\Sigma(\al,\lambda\beta)=\Sigma(\al,\beta),\, \forall \lambda>0,\qquad
(ii)\,\Sigma(\al,-\beta)=-\Sigma(\al,\beta),
$$
and
$$
(iii)\,\sup_\al|\Sigma(\al,\beta)|\leq 8 \|\grad X\|^3_{L^\infty}\||N|^{-1}\|^3_{L^\infty},
$$
as a consequence of \eqref{sioe}.

 Here we have a singular integral operator with odd kernel (see \cite{DP} and \cite{St3}) and therefore a bounded linear map on $L^2(\R^2)$ giving us
$$
\|J_6\|_{L^2}\leq C\|\grad X\|^4_{L^\infty}\||N|^{-1}\|^3_{L^\infty}
\|\partial_{\alpha_j}f\|_{L^2}.
$$

For the family of operators $\mathcal{T}_2(f)(\al)$ we use the splitting $\mathcal{T}_2(f)=I_3+I_4$ where
$$
I_3=\int_{|\beta|>1}\frac{\partial_{\alpha_j}X_k(\al)-\partial_{\alpha_j}X_k(\al-\beta)}{|X(\al)-X(\al-\beta)|^3}f(\al-\beta)d\beta.
$$
Easily we get
$$
I_3\leq 2\|X-(\al,0)\|_{C^{1}}\|F(X)\|^3_{L^\infty} \int_{|\beta|>1}\frac{|f(\al-\beta)|}{|\beta|^3}d\beta,
$$
 while for $I_4$ we proceed with the same method used with $I_2$,  replacing now $X_k(\al)$ by $\partial_{\alpha_j}X_k(\al)$ and $\partial_{\alpha_j}f(\al-\beta)$ by $f(\al-\beta)$.

Next we shall show that the operator $\mathcal{T}_3$ behaves like $\Lambda = (-\Delta)^{\frac12}$. To do that we split it as $I_5 + I_6$ where
$$
I_5=\int_{|\beta|>1}\frac{f(\al)-f(\al-\beta)}{|X(\al)-X(\al-\beta)|^3}d\beta,
$$
 can be easily estimated by
$$
I_5\leq \|F(X)\|^3_{L^\infty}(2\pi|f(\al)|+\int_{|\beta|>1}\frac{|f(\al-\beta)|}{|\beta|^3}d\beta).
$$
The other term is written in the form  $I_6=J_7+J_8$ where
$$
J_7=\int_{|\beta|<1}[\frac{1}{|X(\al)-X(\al-\beta)|^3}-\frac{1}{|\grad X(\al)\cdot \beta|^3}](f(\al)-f(\al-\beta))d\beta.
$$
The identity
$$
f(\al)-f(\al-\beta)=\beta\cdot \int_0^1 \grad f(\al+(s-1)\beta)ds
$$
allows us to treat $J_7$ as we did with $J_5$. To estimate $J_8$  the equality
\begin{equation}\label{lo}
\frac{1}{|\grad X(\al)\cdot\beta|^3}=-\dbu\big(\frac{\beta_1}{|\grad X(\al)\cdot\beta|^3}\big)-\dbd\big(\frac{\beta_2}{|\grad X(\al)\cdot\beta|^3}\big)
\end{equation}
will be very useful. After a careful integration by parts it yields
$$
J_8=PV\int_{|\beta|<1}\frac{\grad f(\al-\beta)\cdot\beta}{|\grad X(\al)\cdot \beta|^3}d\beta-\int_{|\beta|=1}\frac{(f(\al)-f(\al-\beta))|\beta|}{|\grad X(\al)\cdot \beta|^3}dl(\beta).
$$
The principal value in $J_8$ is treated with the same method used for   $J_6$ and since the integral on the circle is inoffensive, so long as $|N|^{-1}$ is in $L^\infty$, the estimate for $\mathcal{T}_3$ follows.

 For the remaining operator one  integrates by parts to get $\mathcal{T}_4=I_7+I_8$ where
$$
I_7=PV\int_{\R^2}P_1(\al,\beta)f(\al-\beta)d\beta,\qquad I_8=PV\int_{\R^2}P_2(\al,\beta)f(\al-\beta)d\beta
$$
with
$$
P_1(\al,\beta)=\frac{\daj X(\al)}{|\grad X(\al)\cdot\beta|^3}-\frac{\daj X(\al-\beta)}{|X(\al)-X(\al-\beta)|^3}
$$
and
\begin{align*}
P_2(\al,\beta)&=
3\frac{(X(\al)-X(\al-\beta)) (X(\al)-X(\al-\beta))\cdot\daj X(\al-\beta)}{|X(\al)-X(\al-\beta)|^5}\\
&\quad-3\frac{\grad X(\al)\cdot\beta ((\grad X(\al)\cdot\beta)\cdot \daj X(\al))}{|\grad X(\al)\cdot\beta|^5}.
\end{align*}
Next we will show how to treat  $I_7$, because the estimate for $I_8$ follows similarly. In $P_1$ we introduce the  further decomposition: $P_1=Q_1+Q_2$ where
$$
Q_1=\daj X(\al)[\frac{1}{|\grad X(\al)\cdot\beta|^3}-\frac{1}{|X(\al)-X(\al-\beta)|^3}],\quad Q_2=\frac{\daj X(\al)-\daj X(\al-\beta)}{|X(\al)-X(\al-\beta)|^3}.
$$
And since the kernel $Q_2$ has already appeared in the operator $\mathcal{T}_1$,   it only remains to control $J_9$ which is given by
$$
J_9=\daj X(\al) PV\int_{\R^2}Q_1(\al,\beta)f(\al-\beta)d\beta.
$$
The following decomposition
$$J_9=\daj X(\al)\int_{|\beta|>1}d\beta+\daj X(\al)PV\int_{|\beta|<1}d\beta=K_1+K_2$$
shows that the term $K_1$ trivializes. Regarding $K_2$  let us write
$$
Q_1=\frac{(|A|^4+|B|^2|A|^2+|B|^4)(A+B)\cdot(A-B)}{|A|^3|B|^3(|A|^3+|B|^3)}
$$
where
$$
A(\al,\beta)=X(\al)-X(\al-\beta),\qquad B(\al,\beta)=\grad X(\al)\cdot\beta.
$$
This formula shows that inside $Q_1$ lies a kernel of degree  $-2$. Then let us take $Q_1=S_1+S_2$ where
$$
S_2=\frac{3|B|^4 B\cdot(A-B)}{|B|^9}=\frac{3B\cdot(A-B)}{|B|^5}.
$$
Next we check that the kernel $S_1$ has degree $-1$, and therefore is easy to handle. Finally we have to consider the kernel $S_2$ appearing in the integral $L$
$$
L=3\daj X(\al) PV\int_{|\beta|<1}\frac{(\grad X(\al)\cdot\beta)\cdot(X(\al)-X(\al-\beta)-\grad X(\al)\cdot\beta)}{|\grad X(\al)\cdot\beta|^5}f(\al-\beta)d\beta.
$$
To do that we introduce a further decomposition $L=M_1+M_2$, with
$$
M_1\!=\!3\daj X(\al)\!\int_{|\beta|<1}\!\!\!\!\!\!\!\!\!\frac{(\grad X(\al)\!\cdot\!\beta)\!\cdot\!(X(\al)\!-\!X(\al\!-\!\beta)\!-\!\grad X(\al)\!\cdot\!\beta\!-\!\frac12 \beta\!\cdot\!\grad^2X(\al)\!\cdot\!\beta)}{|\grad X(\al)\cdot\beta|^5}f(\al\!-\!\beta)d\beta
$$
and
$$
M_2=\frac32\daj X(\al)PV\int_{|\beta|<1}\frac{(\grad X(\al)\cdot\beta)\cdot( \beta\cdot\grad^2X(\al)\cdot\beta)}{|\grad X(\al)\cdot\beta|^5}f(\al-\beta)d\beta,
$$
where $\frac12 \beta\cdot\grad^2X(\al)\cdot\beta$ is the second order term in the Taylor expansion of $X$. It is now easy to check that
$$
M_1\leq C\|\grad X\|^5_{L^\infty}\|X-(\al,0)\|_{C^{2,\delta}}\||N|^{-1}\|^4_{L^\infty}\int_{|\beta|<1}
\frac{|f(\al\!-\!\beta)|}{|\beta|^{2-\delta}}d\beta.
$$
Then we also check that  $M_2$ is controlled like  $J_6$ throughout  the estimate
$$
\|M_2\|_{L^2}\leq C\|\grad X\|^5_{L^\infty}\|\grad^2X\|_{L^{\infty}}\||N|^{-1}\|^4_{L^\infty}\|f\|_{L^2}
$$
which allows us to finish the proof.

\begin{rem}
Having obtained the a priori bounds of the precedent sections, we are in position to implement successfully the same approximation scheme developed in \cite{ADP} to conclude local existence.
\end{rem}

\subsection*{{\bf Acknowledgments}}

\smallskip

We are glad to thank C. Kenig for several wise comments, which help us to simplify our original proof of lemma 5.3.

AC was partially supported by {\sc MTM2008-038} project of the MCINN (Spain).
DC and FG were partially supported by {\sc MTM2008-03754} project of the MCINN (Spain) and StG-203138CDSIF grant of the ERC.
FG was partially supported by NSF-DMS grant 0901810.

\begin{quote}
\begin{tabular}{l}
\textbf{Antonio C\'ordoba} \\
{\small Instituto de Ciencias Matem\'aticas-CSIC-UAM-UC3M-UCM}\\
{\small \& Departamento de Matem\'aticas}\\{\small Facultad de
Ciencias} \\ {\small Universidad Aut\'onoma de Madrid}
\\ {\small Crta. Colmenar Viejo km.~15,  28049 Madrid,
Spain} \\ {\small Email: antonio.cordoba@uam.es}
\end{tabular}
\end{quote}
\begin{quote}
\begin{tabular}{ll}
\textbf{Diego C\'ordoba} &  \textbf{Francisco Gancedo}\\
{\small Instituto de Ciencias Matem\'aticas} & {\small Department of Mathematics}\\
{\small Consejo Superior de Investigaciones Cient\'ificas} & {\small University of Chicago}\\
{\small Serrano 123, 28006 Madrid, Spain} & {\small 5734 University Avenue, Chicago, IL 60637}\\
{\small Email: dcg@icmat.es} & {\small Email: fgancedo@math.uchicago.edu}
\end{tabular}
\end{quote}

\begin{thebibliography}{99}

\bibitem{Ambrose} D. Ambrose, Well-posedness of two-phase Darcy flow in 3D.
\emph{Quart. Appl. Math.}  65,  no. 1, 189-203, 2007.

\bibitem{bear} J. Bear, Dynamics of Fluids in Porous Media, \emph{American
Elsevier}, New York, 1972.

\bibitem{Peter} P. Constantin and M. Pugh. Global solutions for small data to the
Hele-Shaw problem. \emph{Nonlinearity}, 6, 393-415, 1993.

\bibitem{CC}  A. C\'{o}rdoba and D. C\'{o}rdoba, A pointwise estimate for
fractionary derivatives with applications to P.D.E., \textit{Proc.
Natl. Acad. Sci., } 100, 26, 15316-15317, 2003.

\bibitem{ADP} A. C\'ordoba, D. C\'ordoba and F. Gancedo. Interface evolution: the Hele-Shaw and Muskat problems. To appear in \emph{Annals of Math}, 2010.

\bibitem{ADP2} A. C\'ordoba, D. C\'ordoba and F. Gancedo.  The Rayleigh-Taylor condition for the evolution of irrotational fluid interfaces. \emph{Proc. Natl. Acad. Sci.}, 106, no. 27, 10955-10959, 2009.

\bibitem{ADP3} A. C\'ordoba, D. C\'ordoba and F. Gancedo. On the uniqueness for SQG patches, Muskat and Water waves.
Preprint.

\bibitem{DP} D. C\'ordoba and F. Gancedo. Contour dynamics of incompressible 3-D fluids in a porous medium with different densities. \emph{Comm. Math. Phys.} 273, 2, 445-471, 2007.

\bibitem{DP2} D. C\'ordoba and F. Gancedo. A maximum principle for the Muskat problem for fluids with different densities. \emph{Comm. Math. Phys.}, 286, 2, 681-696, 2009.

\bibitem{Darcy} H. Darcy. Les Fontaines Publiques de la Ville de Dijon. \newblock\emph{Dalmont}, Paris, 1856.

\bibitem{ES} J. Escher and G. Simonett. Classical solutions for Hele-Shaw models with surface tension. \emph{Adv. Differential Equations}, 2:619-642, 1997.

\bibitem{lewy} H. Lewy. On the boundary behavior of minimal surfaces. \emph{Proc. Nat. Acad. Sci. U. S. A.} 37, (1951). 103--110.

\bibitem{Muskat} M. Muskat. The flow of homogeneous fluids through porous media. \newblock \emph{New York}, 1937.

\bibitem{S-T} P.G. Saffman and Taylor.
\newblock The penetration of a fluid into a porous medium or Hele-Shaw cell containing a more viscous liquid.
\newblock \emph{Proc. R. Soc. London, Ser. A} 245, 312-329, 1958.


\bibitem{Sanchez} E. Sanchez-Palencia. Homogenization techniques for composite media. Papers from the course held in Udine, July 1--5, 1985. Edited by E. Sánchez-Palencia and A. Zaoui.  \emph{Lecture Notes in Physics}, 272. Springer-Verlag, Berlin, 1987.


\bibitem{SCH} M. Siegel, R. Caflisch and S. Howison. Global
Existence, Singular Solutions, and Ill-Posedness for the Muskat
Problem. \emph{Comm. Pure and Appl. Math.}, 57: 1374-1411, 2004.



\bibitem{St3} E. Stein. Harmonic Analysis. \emph{
Princeton University Press.} Princeton, NJ, 1993.

\bibitem{Tartar} L. Tartar. Incompressible fluid flow in a porous medium -
Convergente of the homogenization processs. Apendix in the book \emph{Nonhomogeneous media and vibration theory} by E. S\'anchez-Palencia. \emph{Lecture Notes in Physics}, 127. Springer-Verlag, Berlin-New York, 1980.

\end{thebibliography}
\end{document}